# Skeleton-stabilized IsoGeometric Analysis: High-regularity Interior-Penalty methods for incompressible viscous flow problems

Tuong Hoang[a,b,*], Clemens V. Verhoosel[a], Ferdinando Auricchio[c],
E. Harald van Brummelen[a], Alessandro Reali[c,d]

[a]*Eindhoven University of Technology – Department of Mechanical Engineering,
P.O. Box 513, 5600MB Eindhoven, The Netherlands*

[b]*IUSS – Istituto Universitario di Studi Superiori Pavia, 27100 Pavia, Italy*

[c]*University of Pavia – Department of Civil Engineering and Architecture, 27100 Pavia, Italy*

[d]*Technische Universität München – Institute for Advanced Study, 85748 Garching, Germany*

## Abstract

A Skeleton-stabilized IsoGeometric Analysis (SIGA) technique is proposed for incompressible viscous flow problems with moderate Reynolds number. The proposed method allows utilizing identical finite dimensional spaces (with arbitrary B-splines/NURBS order and regularity) for the approximation of the pressure and velocity components. The key idea is to stabilize the jumps of high-order derivatives of variables over the skeleton of the mesh. For B-splines/NURBS basis functions of degree $k$ with $C^\alpha$-regularity ($0 \leq \alpha < k$), only the derivative of order $\alpha + 1$ has to be controlled. This stabilization technique thus can be viewed as a high-regularity generalization of the (Continuous) Interior-Penalty Finite Element Method. Numerical experiments are performed for the Stokes and Navier-Stokes equations in two and three dimensions. Oscillation-free solutions and optimal convergence rates are obtained. In terms of the sparsity pattern of the algebraic system, we demonstrate that the block matrix associated with the stabilization term has a considerably smaller bandwidth when using B-splines than when using Lagrange basis functions, even in the case of $C^0$-continuity. This important property makes the proposed isogeometric framework practical from a computational effort point of view.

## 1. Introduction

Isogeometric analysis (IGA) was introduced by Hughes *et al.* [1] as a novel analysis paradigm targeting better integration of Computer Aided Design (CAD) and Finite Element Analysis (FEA). The pivotal idea of IGA is that it directly inherits its basis functions from CAD modeling, where Non-uniform Rational B-splines (NURBS) are the industry standard. For analysis-suitable CAD models, geometrically exact analyses can be performed on the coarsest level of the CAD geometry. This contrasts with conventional FEA, which typically uses Lagrange polynomials as basis functions defined on a geometrically approximate mesh. An additional highly appraised property of IGA is that splines allow one to achieve higher-order continuity, in contrast to the $C^0$-continuity of conventional FEA. We refer to [2, 3] for an overview of established IGA developments.

*Corresponding author.
Email addresses:
tuong.hoang@iusspavia.it, t.hoang@tue.nl (T. Hoang)
c.v.verhoosel@tue.nl (C.V. Verhoosel),
auricchio@unipv.it (F. Auricchio),
e.h.v.brummelen@tue.nl (E.H. van Brummelen),
alessandro.reali@unipv.it (A. Reali)

In the context of viscous flow problems – particularly in the incompressible regime – IGA has been applied very successfully. Within the framework of inf-sup stable spaces for mixed formulations [4], a variety of compatible discretizations has been developed, most notably: Taylor-Hood elements [5, 6, 7], Nédélec elements [6], subgrid elements [8, 7], and H(div)-conforming elements [6, 9, 10, 11]. The mixed discretization approach leads to a saddle point system where the discrete velocity and pressure spaces are chosen differently in order to satisfy the discrete inf-sup condition. The advantage of this approach is that a stable discrete system is obtained straightforwardly from the continuous weak formulation (without any modifications) if the pair of discrete spaces is chosen appropriately.

In practice, employing the same discrete space for the velocity and pressure fields can provide advantages in terms of implementation and computer resources. These advantages become more pronounced in multi-physics problems with many different field variables, for which the derivation of inf-sup stable discrete spaces can be non-trivial. The data structures required to represent the different spaces can make this approach impractical in terms of implementation and computational expenses. Moreover, in the context of IGA, using the same discretization space for all field variables enables direct usage of the CAD basis functions, which is highly beneficial from the vantage point of CAD/FEA integration.

Although there are merits in using the same discrete space for all field variables, without modification this generally leads to an unstable system in the Babuška-Brezzi sense. A common remedy to circumvent this issue is to use stabilization techniques. Various stabilization techniques have been studied in the IGA setting, most notably: Galerkin-least squares and Douglas-Wang stabilization [5] and variational multiscale stabilization (VMS) [12]. The structure of these approaches is that the stabilization is based on element-by-element residuals. We note that recently a combination of VMS and compatible B-splines is studied in [13]. It is also noteworthy that for incompressible elasticity the use of inf-sup stable discretizations can be circumvented by using stream functions [14], the $B$-bar method [15] and the $B^D$-bar method [16].

In this contribution we propose a novel skeleton-based stabilization technique for isogeometric analysis of viscous flow problems, like those described by the Stokes equations and incompressible Navier-Stokes equations with moderate Reynolds numbers. The skeleton-based stabilization allows utilizing identical finite dimensional spaces for the approximation of the pressure and velocity fields. The central idea is to supplement the variational formulation with a consistent penalization term for the jumps of high-order derivatives of the pressure across element interfaces. By taking into account the local continuity at each element interface, the stabilized formulation can be applied to B-splines/NURBS with varying regularities, including the case of multi-patch geometries.

The proposed stabilization technique only controls the $(\alpha + 1)$-th order derivative in the case of B-splines/NURBS basis functions of degree $k$ with $C^\alpha$-regularity. On the one hand, the proposed stabilization technique can be regarded as a generalization of the continuous interior penalty finite element method [17] where $C^0$ Lagrange basis functions are employed. On the other hand, under the minimal stabilization framework [18], we can interpret that the proposed method is related to inf-sup stable approaches investigated in Ref. [7]. This new technique enables the consideration of a large class of problems in isogeometric analysis for fluid flows. The present work encompasses a detailed study of the effect of the stabilization operator on the sparsity pattern of the mixed matrix – including an analysis of its complexity with respect to the B-splines/NURBS order – from which it is observed that the proposed technique optimally exploits the higher-order continuity properties of isogeometric analysis. We present a series of detailed numerical benchmark simulations to demonstrate the effectivity of the stabilization technique. In particular we show that oscillation-free solutions are attained, and the method yields optimal convergence rates under mesh refinements.

The outline of the paper is as follows. In Section 2 we recall the essential aspects of isogeometric analysis. In particular we introduce the skeleton structure and jump operators, and we discuss the local continuity properties across element interfaces. The skeleton-based isogeometric analysis technique for the Navier-Stokes equations is then introduced in Section 3. In Section 4 we discuss the matrix form and implementation aspects of the method, along with a study of the effect of the skeleton-stabilization operator on the sparsity pattern of the algebraic system. A series of numerical test cases is considered in Section 5 to demonstrate the performance of the proposed method. Conclusions are finally presented in Section 6.

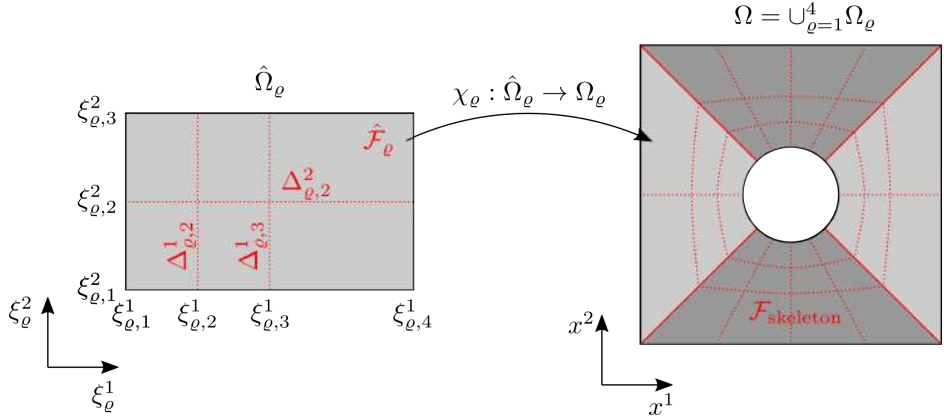

Figure 1: Notations for a parameterization of a multipatch geometry

## 2. Fundamentals of skeleton-based isogeometric analysis

To provide a setting for the skeleton-based stabilization proposed in Section 3 and to introduce the main notational conventions, we first present multi-patch non-uniform rational B-spline (NURBS) spaces. We consider a domain $\Omega \subset \mathbb{R}^d$ (with $d = 2$ or 3) with Lipschitz boundary $\partial \Omega$ as exemplified in Figure 1. The domain $\Omega$ is parameterized by a, possibly multi-patch ($n_{\text{patch}} \geq 1$), non-uniform rational B-spline (NURBS) such that

$$\Omega = \bigcup_{\varrho=1}^{n_{\text{patch}}} \chi_\varrho \circ \widehat{\Omega}_\varrho, \tag{1}$$

where $\widehat{\Omega}_\varrho$ and $\chi_\varrho$ are the patch-wise geometric maps and parameter domains, respectively, with the parametric map defined as

$$\begin{cases} \chi_\varrho : \widehat{\Omega}_\varrho \to \Omega_\varrho, \\ \mathbf{x} = \sum_{I=1}^{n_\varrho} \widehat{R}_{\varrho,I}(\boldsymbol{\xi}_\varrho) \mathbf{X}_{\varrho,I}, \end{cases} \tag{2}$$

where $\{\widehat{R}_{\varrho,I} : \widehat{\Omega}_\varrho \to \mathbb{R}\}_{I=1}^{n_\varrho}$ and $\{\mathbf{X}_{\varrho,I} \in \mathbb{R}^d\}_{I=1}^{n_\varrho}$ are the set of NURBS basis functions and the associated set of control points, respectively. The NURBS basis functions are constructed based on a set of non-decreasing knot vectors, $\{\Xi_\varrho^\delta\}_{\delta=1}^d$, with

$$\Xi_\varrho^\delta = [\underbrace{\xi_{\varrho,1}^\delta, \ldots, \xi_{\varrho,1}^\delta}_{r_{\varrho,1}^\delta \text{ times}}, \underbrace{\xi_{\varrho,2}^\delta, \ldots, \xi_{\varrho,2}^\delta}_{r_{\varrho,2}^\delta \text{ times}}, \ldots, \underbrace{\xi_{\varrho,m_\varrho^\delta}^\delta, \ldots, \xi_{\varrho,m_\varrho^\delta}^\delta}_{r_{\varrho,m_\varrho^\delta}^\delta \text{ times}}], \tag{3}$$

such that the number of basis functions per patch is $n_\varrho = \otimes_{\delta=1}^d \{(\sum_{i=1}^{m_\varrho^\delta} r_{\varrho,i}^\delta) - k_\varrho^\delta - 1\}$, with $k_\varrho^\delta$ the degree of the spline in the direction $\delta$ ($\delta = 1, \ldots, d$). Note that for open B-splines the multiplicity of the first and last knot values is equal to $r_{\varrho,1}^\delta = r_{\varrho,m_\varrho^\delta}^\delta = k_\varrho^\delta + 1$. The regularity of the basis in the parametric directions depends on the order and the multiplicity of the knot value: $\alpha_{\varrho,i}^\delta = k_\varrho^\delta - r_{\varrho,i}^\delta$ for $i = 1, \ldots, m_\varrho^\delta$. On every patch the knot vectors partition the domain into a parametric mesh $\widehat{\mathcal{T}}_\varrho$. The corresponding partitioning of the domain $\Omega$ follows as

$$\mathcal{T}^h = \bigcup_{\varrho=1}^{n_{\text{patch}}} \chi_\varrho \circ \widehat{\mathcal{T}}_\varrho. \tag{4}$$

The superscript $h$ indicates the dependence of the partition on a mesh (resolution) parameter $h > 0$. We associate with the mesh $\mathcal{T}^h$ the skeleton[1]:

$$\mathcal{F}^h_{\text{skeleton}} = \{\partial K \cap \partial K' \mid K, K' \in \mathcal{T}^h, K \neq K'\}. \tag{5}$$

Note that since the skeleton-based stabilization technique considered in this work pertains to inter-element continuity properties, the boundary faces are not incorporated in the skeleton. The skeleton (5) can be decomposed in the intra-patch skeleton, $\mathcal{F}^h_{\text{intra}}$, and the inter-patch skeleton, $\mathcal{F}^h_{\text{inter}}$:

$$\mathcal{F}^h_{\text{intra}} := \bigcup_{\varrho=1}^{n_{\text{patches}}} \chi_\varrho \circ \widehat{F}_\varrho \quad \text{with} \quad \widehat{F}_\varrho := \left\{\partial \widehat{K} \cap \partial \widehat{K'} \mid \widehat{K}, \widehat{K'} \in \widehat{\mathcal{T}}_\varrho, \ \widehat{K} \neq \widehat{K'}\right\}, \tag{6a}$$

$$\mathcal{F}^h_{\text{inter}} := \mathcal{F}^h_{\text{skeleton}} \setminus \mathcal{F}^h_{\text{intra}}. \tag{6b}$$

It evidently follows from these definitions that $\mathcal{F}^h_{\text{skeleton}} = \mathcal{F}^h_{\text{intra}} \cup \mathcal{F}^h_{\text{inter}}$ and $\mathcal{F}^h_{\text{intra}} \cap \mathcal{F}^h_{\text{inter}} = \emptyset$.

Continuity across a patch interface is achieved by matching the knot vectors associated with the two sides of the interface, and by making the corresponding control points on both patches coincident. In terms of the NURBS basis this is equivalent to linking the NURBS basis functions corresponding to the coincident control points. We denote the set of all basis functions over the domain $\Omega$ – where interface functions have been linked – by $\mathcal{R} := \{R_I : \Omega \to \mathbb{R}\}_{I=1}^n$. The space spanned by this basis is denoted by $\mathcal{S} := \text{span}(\mathcal{R})$. Let us note that in the general case of a non-conforming muti-patch structure, multi-patch coupling techniques can be used such as the Nitsche's method [19, 20] or the isogeometric mortar method [21].

To define the regularity of the spline space $\mathcal{S}$ we introduce the plane (or line in the two-dimensional case) in the parameter domain of patch $\varrho$ which is perpendicular to the $\delta$-direction, with its coordinate $\xi^\delta$ equal to that of the knot value $\xi^\delta_{\varrho,i}$ (see Figure 1):

$$\Delta^\delta_{\varrho,i} := \left\{\boldsymbol{\xi} = (\xi^1, \ldots, \xi^d) \mid \xi^\delta = \xi^\delta_{\varrho,i} \text{ and } \xi^{\delta'} \in [\xi^{\delta'}_{\varrho,1}, \xi^{\delta'}_{\varrho,m_\varrho^\delta}] \text{ for } \delta' \neq \delta\right\}. \tag{7}$$

The regularity of the space $\mathcal{S}$ across an intra-patch face $F \in \mathcal{F}^h_{\text{intra}}$ can then be defined through the unique combination of the patch index $\varrho$, the direction $\delta$, and the knot index $i$, such that the associated parametric face $\widehat{F}_\varrho$ resides in the plane $\Delta^\delta_{\varrho,i}$. In combination with the $C^0$-continuity condition across patch boundaries, the regularity of the faces $F \in \mathcal{F}^h_{\text{skeleton}}$ is then given by:

$$\alpha(F) := \begin{cases} \alpha^\delta_{\varrho,i}, & \exists!(\varrho, \delta, i) : \chi_\varrho^{-1} \circ F \subset \Delta^\delta_{\varrho,i}, \quad F \in \mathcal{F}^h_{\text{intra}}, \\ 0, & F \in \mathcal{F}^h_{\text{inter}}. \end{cases} \tag{8}$$

For all functions $f \in \mathcal{S}$ the jumps of its $k$-th normal derivatives across an interface vanish in accordance with

$$[\![\partial_n^k f]\!] = 0, \ 0 \leq k \leq \alpha(F), \tag{9}$$

where the jump for some function $\phi$ is defined as $[\![\phi]\!] \equiv [\![\phi]\!]_F := \phi^+ - \phi^-$, and the superscripts $+$ and $-$ refer to the traces of $\phi$ on the two opposite sides of $F$.

From (8) it is inferred that in the interior of a patch the regularity per direction is controlled by the knot vector multiplicity, while across patch boundaries merely $C^0$-continuity of the basis holds. We denote by $\mathcal{S}^k_{h,\alpha} \equiv \mathcal{S}^k_\alpha$ the spline space with mesh size index $h$, global isotropic degree $k$ and per skeleton face regularity $\alpha$ in accordance with definition (8). In the special case of a global intra-patch regularity $\bar{\alpha} \in \mathbb{N}$, i.e., $\alpha(F) = \bar{\alpha}$, $0 \leq \bar{\alpha} \leq k-1$ $\forall F \in \mathcal{F}^h_{\text{intra}}$ we denote the function space by $\mathcal{S}^k_{\bar{\alpha}}$. A special case of this function space is that in which full regularity is achieved, i.e., $\bar{\alpha} = k-1$.

---

[1]This should not be confused with the topological skeleton concept in geometric modeling.

## 3. Skeleton-stabilized Isogeometric Analysis for the Navier-Stokes equations

In this section we introduce the skeleton-penalty formulation for the Navier-Stokes equations in the context of Isogeometric Analysis. We commence with the formulation of the time-dependent Navier-Stokes equations in Section 3.1. Next, we introduce the discrete skeleton-penalty formulation in Section 3.2.

### 3.1. The time-dependent Navier-Stokes equations

We consider the unsteady incompressible Navier-Stokes equations on the open bounded domain $\Omega \in \mathbb{R}^d$ (with $d = 2$ or $3$). The Lipschitz boundary $\partial\Omega$ is split in two complementary open subsets $\Gamma_D$ and $\Gamma_N$ (such that $\overline{\Gamma_D} \cup \overline{\Gamma_N} = \partial\Omega$ and $\Gamma_D \cap \Gamma_N = \emptyset$) for Dirichlet and Neumann conditions, respectively. The outward-pointing unit normal vector to $\partial\Omega$ is denoted by $\mathbf{n}$. For any time instant $t \in [0,T]$ the Navier-Stokes equations for the velocity field $\mathbf{u} : \Omega \times [0,T] \to \mathbb{R}^d$ and pressure field $p : \Omega \times [0,T] \to \mathbb{R}$ read:

$$\begin{cases} \text{Find } \mathbf{u} : \Omega \times [0,T] \to \mathbb{R}^d, \text{ and } p : \Omega \times [0,T] \to \mathbb{R} \text{ such that:} \\ \partial_t \mathbf{u} + \nabla \cdot (\mathbf{u} \otimes \mathbf{u}) - \nabla \cdot (2\mu \nabla^s \mathbf{u}) + \nabla p = \mathbf{f} \quad \text{in } \Omega \times (0,T), \\ \nabla \cdot \mathbf{u} = 0 \quad \text{in } \Omega \times (0,T), \\ \mathbf{u} = \mathbf{0} \quad \text{on } \Gamma_D \times (0,T), \\ 2\mu \nabla^s \mathbf{u} \cdot \mathbf{n} - p\mathbf{n} = \mathbf{h} \quad \text{on } \Gamma_N \times (0,T), \\ \mathbf{u} = \mathbf{u}_0 \quad \text{in } \Omega \times \{0\}. \end{cases} \tag{10}$$

Here $\mu$ represents the kinematic viscosity, and the symmetric gradient of the velocity field is denoted by $\nabla^s \mathbf{u} := \frac{1}{2}\left(\nabla \mathbf{u} + (\nabla \mathbf{u})^T\right)$. The exogenous data $\mathbf{f} : \Omega \times (0, \infty) \to \mathbb{R}^d$ and $\mathbf{h} : \Gamma_N \times (0, \infty) \to \mathbb{R}^d$, represent the body forces and Neumann conditions, respectively. Without loss of generality we herein assume the Dirichlet data to be homogeneous. The initial conditions in (10) are denoted by $\mathbf{u}_0 : \Omega \to \mathbb{R}^d$.

For any vector space $\mathcal{V}$, we denote by $\mathcal{L}(0,T;\mathcal{V})$ a suitable linear space of $\mathcal{V}$-valued functions on the time interval $(0,T)$. We consider the following weak formulation of (10):

$$\begin{cases} \text{Find } \mathbf{u} \in \mathcal{L}(0,T;\boldsymbol{\mathcal{V}}_{0,\Gamma_D}) \text{ and } p \in \mathcal{L}(0,T;\mathcal{Q}), \text{ given } \mathbf{u}(0) = \mathbf{u}_0, \\ \text{such that for almost all } t \in (0,T): \\ (\partial_t \mathbf{u}, \mathbf{w}) + c(\mathbf{u}; \mathbf{u}, \mathbf{w}) + a(\mathbf{u}, \mathbf{w}) + b(p, \mathbf{w}) = \ell(\mathbf{w}) \quad \forall \mathbf{w} \in \boldsymbol{\mathcal{V}}_{0,\Gamma_D}, \\ \qquad\qquad\qquad\qquad\qquad\qquad\qquad b(q, \mathbf{u}) = 0 \quad \forall q \in \mathcal{Q}. \end{cases} \tag{11}$$

The trilinear, bilinear, and linear forms in this formulation are defined as

$$c(\mathbf{v}; \mathbf{u}, \mathbf{w}) := (\mathbf{v} \cdot \nabla \mathbf{u}, \mathbf{w}), \tag{12a}$$

$$a(\mathbf{u}, \mathbf{w}) := 2\mu \left(\nabla^s \mathbf{u}, \nabla^s \mathbf{w}\right), \tag{12b}$$

$$b(q, \mathbf{w}) := -(q, \text{div}\mathbf{w}), \tag{12c}$$

$$\ell(\mathbf{w}) := (\mathbf{f}, \mathbf{w}) + \langle \mathbf{h}, \mathbf{w}\rangle_{\Gamma_N}, \tag{12d}$$

where $(\cdot, \cdot)$ and $\langle \cdot, \cdot \rangle_{\Gamma_N}$ denote the inner product in $L^2(\Omega)$ and dual product in $L^2(\Gamma_N)$, respectively. The function spaces in (11) are defined as

$$\boldsymbol{\mathcal{V}}_{0,\Gamma_D} := \left\{\mathbf{u} \in [H^1(\Omega)]^d : \mathbf{u} = \mathbf{0} \text{ on } \Gamma_D\right\}, \qquad \mathcal{Q} := L^2(\Omega).$$

In the case of pure Dirichlet boundary conditions, i.e., if $\Gamma_D$ coincides with all of $\partial\Omega$, the pressure is determined up to a constant. In that case, the pressure space is subject to the zero average pressure condition:

$$\mathcal{Q} := L^2_0(\Omega) \equiv \left\{q \in L^2(\Omega) : \int_\Omega q \, d\Omega = 0\right\}. \tag{13}$$

*3.2. The Isogeometric Skeleton-Penalty method with identical discrete spaces of velocity and pressure*

In this contribution we study the discretization of (11) by utilizing identical spline discretizations for the velocity and pressure fields. The global isotropic order of the spline space is denoted by $k$ and its regularity by $\alpha$ (with $0 \leq \alpha(F) \leq k-1 \; \forall F \in \mathcal{F}^h_{\text{skeleton}}$; see Section 2):

$$\boldsymbol{\mathcal{V}}^h := \left[\mathcal{S}^k_\alpha\right]^d \cap \boldsymbol{\mathcal{V}}_{0,\Gamma_D}, \qquad \mathcal{Q}^h := \mathcal{S}^k_\alpha \cap \mathcal{Q}. \qquad (14)$$

The semi-discretization in space of the weak form (11) then reads:

$$\begin{cases} \text{Find } \mathbf{u}^h \in \mathcal{L}(0,T;\boldsymbol{\mathcal{V}}^h) \text{ and } p^h \in \mathcal{L}(0,T;\mathcal{Q}^h), \text{ given } \mathbf{u}^h(0) = \mathbf{u}^h_0, \\ \text{such that for almost all } t \in (0,T): \\ \begin{aligned} (\partial_t \mathbf{u}^h, \mathbf{w}^h) + c(\mathbf{u}^h;\mathbf{u}^h,\mathbf{w}^h) + a(\mathbf{u}^h,\mathbf{w}^h) + b(p^h,\mathbf{w}^h) &= \ell(\mathbf{w}^h) & \forall \mathbf{w}^h \in \boldsymbol{\mathcal{V}}^h, \\ b(q^h,\mathbf{u}^h) &= 0 & \forall q^h \in \mathcal{Q}^h. \end{aligned} \end{cases} \qquad (15)$$

The pair of spaces $(\boldsymbol{\mathcal{V}}^h, \mathcal{Q}^h)$ in (14) does not satisfy the inf-sup condition, and hence the discretization in (15) is unstable. To stabilize the system, we propose to supplement the formulation with the skeleton-penalty term,

$$s(p^h,q^h) := \sum_{F \in \mathcal{F}^h_{skeleton}} \int_F \gamma \mu^{-1} h_F^{2\alpha+3} [\![\partial_n^{\alpha+1} p^h]\!][\![\partial_n^{\alpha+1} q^h]\!] d\Gamma, \qquad (16)$$

where $\alpha$ is the regularity of the considered spline space at the element interface $F \in \mathcal{F}^h_{skeleton}$, $\gamma > 0$ is a global stabilization parameter, and $h_F$ is a length scale associated with this element interface. Here we define this length scale as

$$h_F := \frac{|K_F^+|_d + |K_F^-|_d}{2|F|_{d-1}}, \qquad (17)$$

where $K_F^+$ and $K_F^-$ are two elements sharing the interface $F$, and $|\cdot|_d$ is the $d$-dimensional Hausdorff measure. The stabilized semi-discrete system – to which we refer as the *isogeometric skeleton-penalty formulation* for the Navier-Stokes equations – then reads:

$$\begin{cases} \text{Find } \mathbf{u}^h \in \mathcal{L}(0,T;\boldsymbol{\mathcal{V}}^h) \text{ and } p^h \in \mathcal{L}(0,T;\mathcal{Q}^h), \text{ given } \mathbf{u}^h(0) = \mathbf{u}^h_0, \\ \text{such that for almost all } t \in (0,T): \\ \begin{aligned} (\partial_t \mathbf{u}^h, \mathbf{w}^h) + c(\mathbf{u}^h;\mathbf{u}^h,\mathbf{w}^h) + a(\mathbf{u}^h,\mathbf{w}^h) + b(p^h,\mathbf{w}^h) &= \ell(\mathbf{w}^h) & \forall \mathbf{w}^h \in \boldsymbol{\mathcal{V}}^h, \\ b(q^h,\mathbf{u}^h) - s(p^h,q^h) &= 0 & \forall q^h \in \mathcal{Q}^h. \end{aligned} \end{cases} \qquad (18)$$

**Remark 1.** *The power $2\alpha+3$ associated with the interface length $h_F$ in (16) follows from scaling arguments. The global stabilization parameter $\gamma$ depends on the utilized spline space $S^p_\alpha$. For a sufficiently smooth pressure solution, viz. $p \in H^{\alpha+1}(\Omega)$, the stabilized formulation (18) is variationally consistent with the weak form (11).*

**Remark 2.** *A special case, which is very common for CAD models, is that in which the highest regularity spline space, $S^k_{k-1}$, is used within each patch of the domain, while $C^0$-continuity is established between patches. The skeleton-penalty term (16) in this case reads:*

$$s(p^h,q^h) := \sum_{F \in \mathcal{F}^h_{intra}} \int_F \gamma \mu^{-1} h_F^{2k+1} [\![\partial_n^k p^h]\!][\![\partial_n^k q^h]\!] d\Gamma + \sum_{F \in \mathcal{F}^h_{inter}} \int_F \gamma \mu^{-1} h_F^3 [\![\partial_n p^h]\!][\![\partial_n q^h]\!] d\Gamma. \qquad (19)$$

**Remark 3.** *The formulation (18) based on the skeleton-penalty stabilization term (16) can also be applied to Lagrange bases, which is – in terms of function spaces – equivalent with the special case corresponding to regularity $\alpha = 0$. In this case, only the jump of first order derivatives must be stabilized. This case is known as the continuous interior penalty finite element method [17]. For higher smoothness B-splines, $S^k_\alpha$, with regularity $\alpha \geq 1$, the jump of first order derivatives vanishes, as a consequence of which the formulation in [17] cannot be applied. Thus, formulation (16) is the high-regularity generalization of the continuous interior penalty finite element method. Note that although the formulation in [17] is equivalent to the special case of $\alpha = 0$, the use of higher-order Bézier elements instead of higher-order Lagrange elements affects the sparsity pattern (see Section 4.3).*

**Remark 4.** *The weak formulation of the steady Stokes problem associated with* (11) *is given by:*

$$\begin{cases} \text{Find } \mathbf{u} \in \boldsymbol{\mathcal{V}}_{0,\Gamma_D} \text{ and } p \in \mathcal{Q} \text{ such that:} \\ a(\mathbf{u},\mathbf{w}) + b(p,\mathbf{w}) = \ell(\mathbf{w}) \quad \forall \mathbf{w} \in \boldsymbol{\mathcal{V}}_{0,\Gamma_D}, \\ \qquad\qquad\quad b(q,\mathbf{u}) = 0 \quad\quad \forall q \in \mathcal{Q}. \end{cases} \tag{20}$$

*Similar to formulation* (18), *the isogeometric skeleton-penalty formulation for the Stokes equations reads:*

$$\begin{cases} \text{Find } \mathbf{u}^h \in \boldsymbol{\mathcal{V}}^h \text{ and } p^h \in \mathcal{Q}^h \text{ such that:} \\ a(\mathbf{u}^h,\mathbf{w}^h) + b(p^h,\mathbf{w}^h) = \ell(\mathbf{w}^h) \quad \forall \mathbf{w}^h \in \boldsymbol{\mathcal{V}}^h, \\ b(q^h,\mathbf{u}^h) - s(p^h,q^h) = 0 \quad\quad \forall q^h \in \mathcal{Q}^h. \end{cases} \tag{21}$$

*It is well-known that problem* (20) *is the first-order optimality condition for the saddle point* $(\mathbf{u},p)$ *of the Lagrangian functional (see e.g. [4])*

$$\mathcal{L}(\mathbf{v},q) = \frac{1}{2}a(\mathbf{v},\mathbf{v}) + b(q,\mathbf{v}) - \ell(\mathbf{v}), \quad (\mathbf{v},q) \in \boldsymbol{\mathcal{V}}_{0,\Gamma_D} \times \mathcal{Q}. \tag{22}$$

*Analogously, the stabilized discrete system* (21) *is related to the optimization problem for the modified Lagrangian functional*

$$\mathcal{L}^h(\mathbf{v}^h,q^h) = \tfrac{1}{2}a(\mathbf{v}^h,\mathbf{v}^h) + b(q,\mathbf{v}^h) - \ell(\mathbf{v}^h) - J(q^h), \quad (\mathbf{v}^h,q^h) \in \boldsymbol{\mathcal{V}}^h \times \mathcal{Q}^h, \tag{23}$$

*with*

$$J(q^h) = \frac{\gamma}{2} \sum_{F \in \mathcal{F}^h_{skeleton}} \int_F \mu^{-1} h_F^{2\alpha+3} \left| [\![\partial_n^{\alpha+1} q^h]\!] \right|^2 \, d\Gamma. \tag{24}$$

*The stabilized discrete system* (21) *follows directly from the first-order optimality condition for this modified Lagrangian functional, and the stabilization term* (16) *appears as the variational derivative of* (24). *From* (24) *it is seen that the stabilization term* (16) *effectively leads to minimization of the jump of high-order derivatives of the pressure over the skeleton* $\mathcal{F}^h_{skeleton}$ *in a least-squares sense.*

**Remark 5.** *To provide a rationale for the proposed skeleton-based stabilization technique, we first note that for* $0 \le \alpha \le k-1$ *the velocity-pressure pair* $(\mathcal{S}^k_{h,\alpha}, \mathcal{S}^k_{h,\alpha+1})$, *in which the regularity of the pressure space exceeds that of the velocity space by 1, is inf-sup stable [7]. The skeleton-based stabilization term* $s(p^h,q^h)$ *in (16) essentially penalizes the deviation of the pressure* $p^h \in \mathcal{S}^k_{h,\alpha}$ *from the stable space* $\mathcal{S}^k_{h,\alpha+1}$. *Indeed, it holds that*

$$s(p^h,p^h) = 0 \qquad \forall p^h \in \mathcal{S}^k_{h,\alpha+1} \subset \mathcal{S}^k_{h,\alpha} \tag{25}$$

$$s(p^h,p^h) > 0 \qquad \forall p^h \in \mathcal{S}^k_{h,\alpha} \setminus \mathcal{S}^k_{h,\alpha+1} \tag{26}$$

*which indicates that inf-sup stability can be restored by adding* $s(\cdot,\cdot)$ *with a properly scaled multiplicative constant to the formulation. The mesh dependence of the stabilization constant according to* $h_F^{2\alpha+3}$ *follows from a simple scaling argument. It is noteworthy that the setting and selection of the stabilization term are in fact reminiscent of minimal stabilizations for mixed problems as presented in [18, §4]. Alternatively, for the maximum-regularity case* $(\alpha = k-1)$, *the stability of the skeleton-stabilized formulation can be related to the inf-sup stability of the maximum-regularity sub-grid element [7, Thm. 4.2] . A proof of inf-sup stability is beyond the scope of this paper. Inf-sup stability of the skeleton-stabilized formulation is investigated numerically in Section 5.1.*

**Remark 6.** *For quasi-uniform meshes, the length scale* $h_F$ *can alternatively be defined as*

$$h_F := \frac{|K_F^+|_d^{1/d} + |K_F^-|_d^{1/d}}{2}, \tag{27}$$

*or, even simpler, as*

$$h_F := \begin{cases} length(F) & d = 2, \\ diam(F) & d = 3. \end{cases} \tag{28}$$

*The numerical results presented in Section 5 are based on definition* (28).

# 4. The algebraic form of Skeleton-stabilized Isogeometric Analysis

In this section we discuss various algorithmic aspects of the proposed skeleton-based stabilization framework. In Section 4.1 we briefly discuss the employed solution procedure for the unsteady Navier-Stokes equations, after which the algebraic form of the formulation is introduced in Section 4.2. The effect of the proposed stabilization term on the sparsity pattern of the system matrix is then studied in detail in Section 4.3.

## 4.1. The unsteady Navier-Stokes solution procedure

We employ a standard solution procedure for the unsteady Navier-Stokes equations. Crank-Nicolson time integration is considered in combination with Picard iterations for solving the nonlinear algebraic problem in each time step. The employed solution strategy is summarized in Algorithm 1. We denote the constant time step size by $\Delta t$ and the time step index by $\imath$, such that $t = \imath \Delta t$. The solution at time step $\imath$ is denoted by $(\mathbf{u}^\imath, p^\imath)$, and the time-dependence of the non-autonomous linear operator $\ell(\mathbf{w})$ is similarly indicated by a superscript: $\ell^\imath(\mathbf{w})$. The Picard iteration counter is denoted by $\jmath$, and the unresolved solution at iteration $\jmath$ by $(\mathbf{u}^\imath_\jmath, p^\imath_\jmath)$. Note that for the sake of notational brevity we here omit the superscript $h$ from the variables.

```
Input: u_0, Δt, tol                              # initial condition, time step, Picard tolerance
# Initialization at t = 0
u^0 = u_0
# Time iteration (θ = 1/2:  Crank-Nicolson)
for ı in 1, 2, ...:
    # Picard iteration
    u_0^ı = u^(ı-1)
    p_0^ı = p^(ı-1) if ı > 1 else 0
    for ȷ in 1, 2, ...:
        ⎧ Find (u_ȷ^ı, p_ȷ^ı) ∈ V^h × Q^h such that:
        ⎪ ((u_ȷ^ı - u^(ı-1))/Δt, w) + θ(c(u_{ȷ-1}^ı; u_ȷ^ı, w) + a(u_ȷ^ı, w))
        ⎨ +(1-θ)(c(u^(ı-1); u^(ı-1), w) + a(u^(ı-1), w)) + b(p_ȷ^ı, w)  =  θℓ^ı(w) + (1-θ)ℓ^(ı-1)(w)   ∀w ∈ V^h,
        ⎩                                                 b(q, u_ȷ^ı) - s(p_ȷ^ı, q)  =  0                              ∀q ∈ Q^h.

        if max{‖u_ȷ^ı - u_{ȷ-1}^ı‖, ‖p_ȷ^ı - p_{ȷ-1}^ı‖} < tol:
          | break
        end
    end
end
```

**Algorithm 1:** Solution procedure for the unsteady Navier-Stokes equations

## 4.2. The algebraic form

Let $\{\mathbf{R}_i\}_{i=1}^{n_u}$ and $\{R_i\}_{i=1}^{n_p}$ denote two sets of NURBS basis functions for the velocity and pressure fields, respectively. The vector-valued velocity basis functions are defined as

$$\mathbf{R}_{i=j+\delta n} = R_j \mathbf{e}_\delta, \qquad j = 1, \ldots, n \text{ and } \delta = 1, \ldots, d \qquad (29)$$

where $n$ is the number of control points, $d$ the number of spatial dimensions (evidently, $n_u = dn$ and $n_p = n$), and $\mathbf{e}_\delta$ is the unit vector in the direction $\delta$. The basis functions span the discrete velocity and pressure spaces

$$\mathcal{V}^h = \text{span}\{\mathbf{R}_i\}_{i=1}^{n_u}, \qquad \mathcal{Q}^h = \text{span}\{R_i\}_{i=1}^{n_p}. \qquad (30)$$

The approximate velocity field $\mathbf{u}^h(\mathbf{x},t)$ and pressure field $p^h(\mathbf{x},t)$ can then be written as

$$\mathbf{u}^h(\mathbf{x},t) = \sum_{i=1}^{n_u} \mathbf{R}_i(\mathbf{x})\hat{u}_i(t), \qquad p^h(\mathbf{x},t) = \sum_{i=1}^{n_p} R_i(\mathbf{x})\hat{p}_i(t), \qquad (31)$$

where $\hat{\mathbf{u}}(t) = (\hat{u}_1, \hat{u}_2, \ldots, \hat{u}_{n_u})^T$ and $\hat{\mathbf{p}}(t) = (\hat{p}_1, \hat{p}_2, \ldots, \hat{p}_{n_p})^T$ are vectors of degrees of freedom. The corresponding algebraic form of (18) then reads

$$\begin{cases} \text{For each } t \in (0,T), \text{ find } \hat{\mathbf{u}} = \hat{\mathbf{u}}(t) \in \mathbb{R}^{n_u} \text{ and } \hat{\mathbf{p}} = \hat{\mathbf{p}}(t) \in \mathbb{R}^{n_p}, \text{ given } \hat{\mathbf{u}}(0) = \hat{\mathbf{u}}_0, \text{ such that:} \\ \mathbf{M}\partial_t \hat{\mathbf{u}} + [\mathbf{C}(\hat{\mathbf{u}}) + \mathbf{A}]\hat{\mathbf{u}} + \mathbf{B}^T \hat{\mathbf{p}} = \mathbf{f}, \\ \mathbf{B}\hat{\mathbf{u}} - \mathbf{S}\hat{\mathbf{p}} = \mathbf{0}. \end{cases} \qquad (32)$$

with the matrix entries given by:

$$A_{ij} = a(\mathbf{R}_j, \mathbf{R}_i), \qquad (33\text{a})$$
$$B_{ij} = b(R_i, \mathbf{R}_j), \qquad (33\text{b})$$
$$C(\hat{\mathbf{u}})_{ij} = c(\hat{\mathbf{u}}; \mathbf{R}_j, \mathbf{R}_i), \qquad (33\text{c})$$
$$S_{ij} = s(R_j, R_i), \qquad (33\text{d})$$
$$M_{ij} = (\mathbf{R}_j, \mathbf{R}_i), \qquad (33\text{e})$$
$$f_i = \ell(\mathbf{R}_i). \qquad (33\text{f})$$

The algebraic form of the solution Algorithm 1 is presented in Algorithm 2.

```
Input: û₀, Δt, tol                           #initial condition vector, time step, Picard tolerance
#Initialization at t = 0
û⁰ = û₀
#Time iteration (θ = ½:  Crank-Nicolson)
for ι in 1, 2, ...:
    #Picard iteration
    û₀ⁱ = ûⁱ⁻¹
    p̂₀ⁱ = p̂ⁱ⁻¹ if ι > 1 else 0
    for ȷ in 1, 2, ...:
        Obtain (ûⱼⁱ, p̂ⱼⁱ) by solving the linear system:

        [ 1/Δt M + θ((C(ûⱼ₋₁ⁱ) + A))   Bᵀ ] [ûⱼⁱ]   [(1/Δt M - (1-θ)(C(ûⁱ⁻¹) + A))ûⁱ⁻¹ + θfⁱ + (1-θ)fⁱ⁻¹]
        [            B                  -S ] [p̂ⱼⁱ] = [                     0                           ]

        if max{‖ûⱼⁱ - ûⱼ₋₁ⁱ‖, ‖p̂ⱼⁱ - p̂ⱼ₋₁ⁱ‖} < tol:
            break
        end
    end
end
```
**Algorithm 2:** Algebraic form of the solution procedure for the unsteady Navier-Stokes equations

We note that computation of the stabilization matrix $\mathbf{S}$ requires a data structure related to the skeleton $\mathcal{F}^h_{skeleton}$ of the mesh $\mathcal{T}^h$. This data structure is constructed such that at each element interface $F \in \mathcal{F}^h_{skeleton}$, the jump of high-order derivatives of the basis functions over $F$ can be evaluated. It should be noted that this skeleton structure is compatible with the recently proposed efficient row-by-row assembly procedure for IGA [22].

### 4.3. The $k/\alpha$-complexity of the skeleton-penalty operator on sparsity pattern

The skeleton-based stabilization operator (16) affects the sparsity pattern of the discretized Navier-Stokes system due to the fact that the jump operators on the (higher-order) derivatives provide additional connectivity between basis functions. To illustrate this effect we consider the spline space

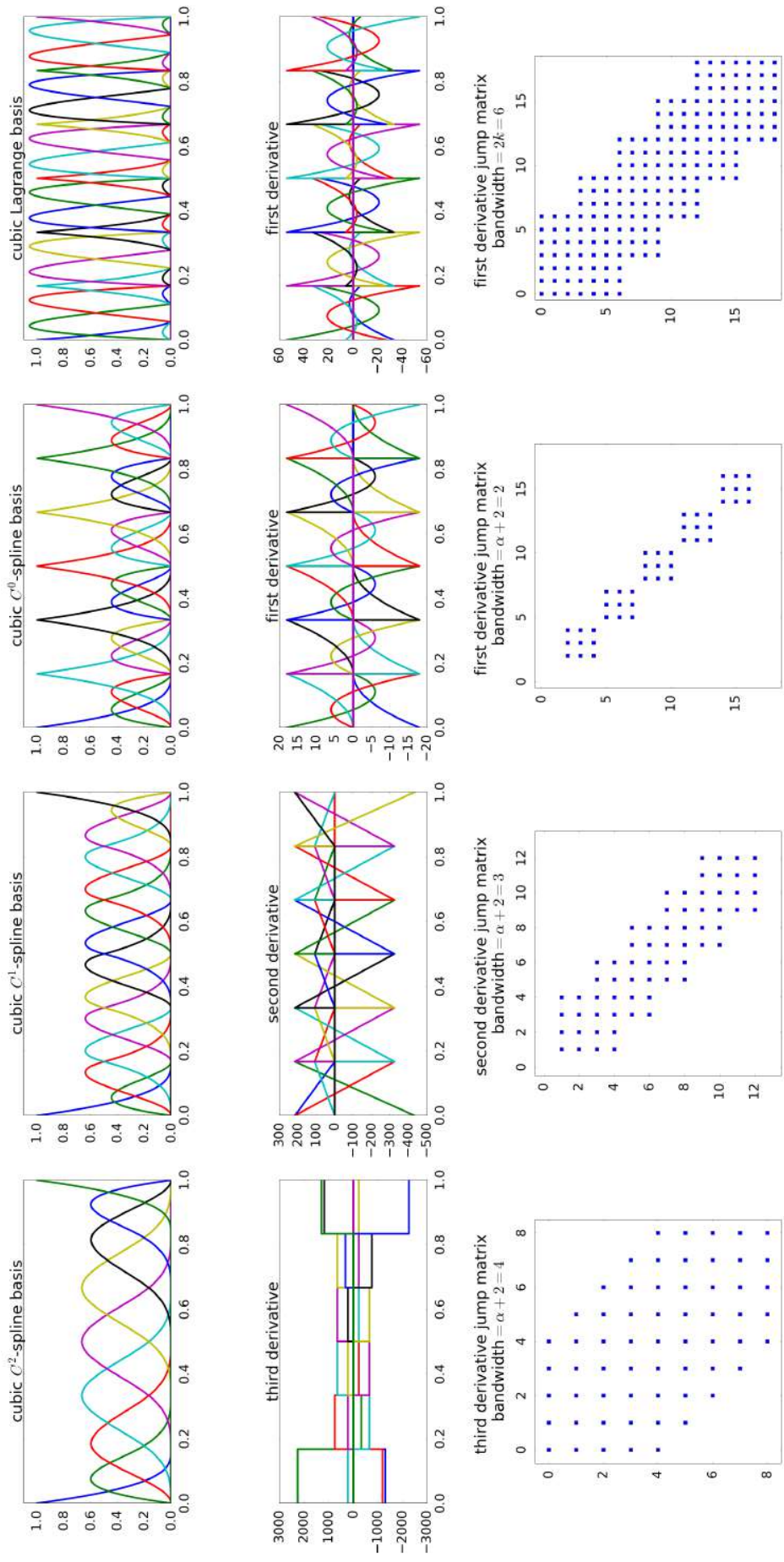

Figure 2: Sparsity pattern of the skeleton-penalty matrix, illustrated with univariate cubic spaces: spline $S_\alpha^3$ space with full regularity $\alpha = 2$ (first column), reduced regularity $\alpha = 1$ (second column), minimal regularity $\alpha = 0$ (third column), and $C^0$ Lagrange space (last column). The top row shows the basis functions, the second row the stabilized (order $\alpha+1$) derivatives, and the third row the matrix sparsity pattern of the skeleton-penalty matrix $\mathbf{S}$. The bandwidths of $\mathbf{S}$ in the spline cases are $\alpha+2$, much smaller than in the Lagrange case $2k$.

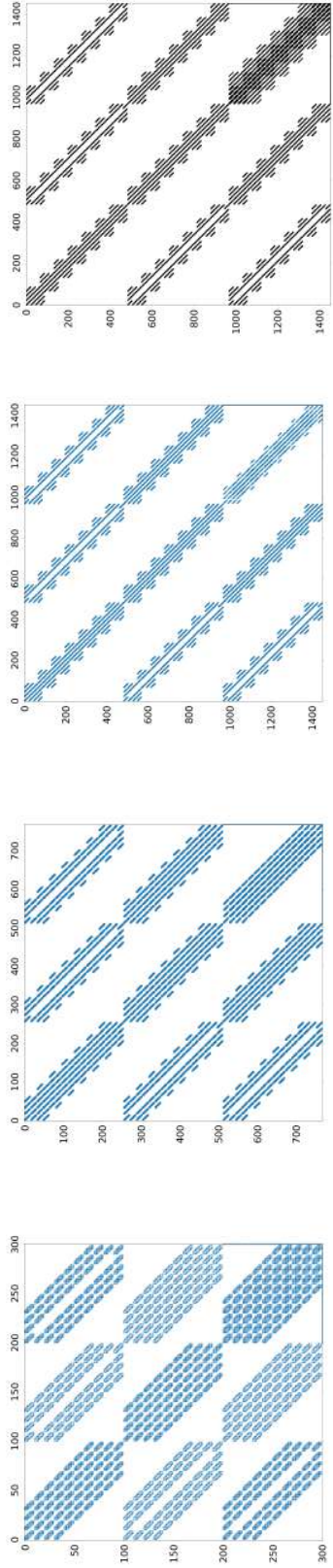
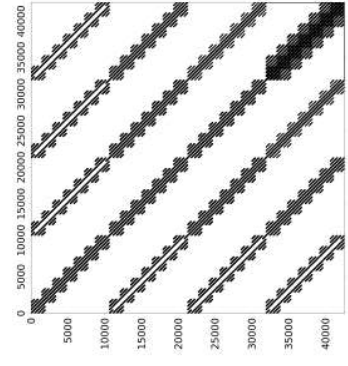
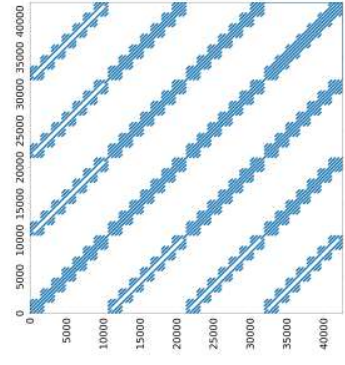
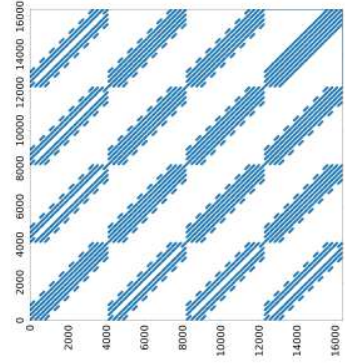
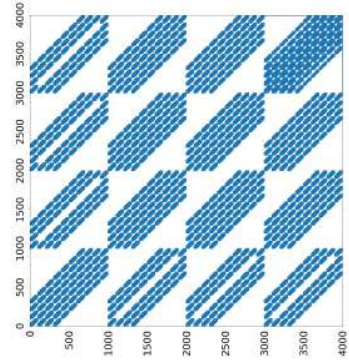

(a) 2D $C^2$ cubic spline
(b) 2D $C^1$ cubic spline
(c) 2D $C^0$ cubic spline
(d) 2D cubic Lagrange
(e) 3D $C^2$ cubic spline
(f) 3D $C^1$ cubic spline
(g) 3D $C^0$ cubic spline
(h) 3D cubic Lagrange

Figure 3: Sparsity pattern of the skeleton-stabilized system in two and three dimensions for various cubic spline spaces $S_\alpha^3$ with: (first column) full regularity $\alpha = 2$, (second column) reduced regularity $\alpha = 1$, (third column) minimal regularity $\alpha = 0$. The fourth column shows the $C^0$ cubic Lagrange space. The bandwidths of the added stabilization block in the spline cases are much smaller than in the Lagrange case, even in the case of $C^0$ regularity. (Note that the size of figures are not scaled with the size of the matrices.)

$\mathcal{S}_\alpha^k$, for which the derivative of order $\alpha + 1$ are stabilized. The top row of Figure 2 shows univariate cubic B-spline bases with $C^2$, $C^1$, $C^0$-regularity, and $C^0$ Lagrange (from left to right). The second row plots the stabilized (order $\alpha + 1$) derivatives for each basis. The third row shows the sparsity pattern of the skeleton-penalty matrix **S** associated with the operator $s(p^h, q^h)$.

The bandwidth[2] of the skeleton-penalty matrix **S** is equal to $\alpha + 2$, which ranges from 2 for $C^0$-splines (or at patch interfaces) to a maximum of $k + 1$ for splines with full continuity (typical for intra patch interfaces). This observed decrease in bandwidth with decrease in regularity stems from the fact that the number of order $\alpha + 1$ derivatives of the basis functions that vanish on the interfaces increases with $\alpha$. This behavior contrasts with classical $C^0$ Lagrange basis functions, for which the bandwidth is equal to $2k$ (the last column of Figure 2). The resulting increase in bandwidth of the jump stabilization matrix with increase in Lagrange basis order is an important drawback of the interior penalty method compared to element-based stabilization techniques. By construction, B-spline bases ameliorate this issue in the sense that even at full continuity the bandwidth of the skeleton-penalty matrix is considerably smaller than that of the Lagrange basis of equal order.

The sparsity patterns of the complete system matrix in two and three dimensions are shown in Figure 3. The skeleton-stabilization term corresponds to the bottom-right block of each matrix. Similar to the observations for the one-dimensional setting, spline bases provide smaller stencils than Lagrange bases. This is an important advantage of the Skeleton-stabiized IGA approach over standard $C^0$ FEM, especially for large systems and in conjunction with iterative solvers.

## 5. Numerical experiments

In this section we investigate the numerical performance of the Skeleton-stabilized IsoGeometric Analysis framework for a range of numerical test cases for viscous flow problems. These test cases focus on various aspects of the framework, most notably its accuracy and convergence under mesh refinement, its stability, and its robustness with respect to the model parameters.

### 5.1. Steady Stokes flow in a unit square

We consider the steady two-dimensional Stokes problem – *i.e.*, problem (11) without time-dependent and convective terms – in the unit square domain $\Omega = (0, 1)^2$. The body force **f** is taken in accordance with the manufactured solution [6]:

$$\mathbf{u} = \begin{pmatrix} 2e^x(-1+x)^2 x^2(y^2-y)(-1+2y) \\ (-e^x(-1+x)x(-2+x(3+x))(-1+y)^2 y^2) \end{pmatrix} \tag{34a}$$

$$p = (-424 + 156e + (y^2 - y)(-456 + e^x(456 + x^2(228 - 5(y^2 - y)) + 2x(-228 + (y^2 - y)) \\ + 2x^3(-36 + (y^2 - y)) + x^4(12 + (y^2 - y))))). \tag{34b}$$

This manufactured solution is visualized in Figure 4a. Note that homogeneous Dirichlet boundary conditions are imposed on the complete boundary $\partial\Omega$, and that a zero average pressure condition, $\int_\Omega p \, d\Omega = 0$, is imposed to establish well-posedness. We use a Lagrange multiplier approach to enforce this condition.

In Figure 4 we study the asymptotic $h$-convergence behavior of the proposed method for B-splines of degree $k = 1, 2, 3$ with the highest possible regularities, *i.e.* $C^{k-1}$. The coarsest mesh considered consists of $4 \times 4$ elements, which is uniformly refined until a $128 \times 128$ mesh is obtained. The stabilization parameter is taken as $\gamma = 1$ ($k = 1$), $5 \times 10^{-2}$ ($k = 2$), $10^{-3}$ ($k = 3$). The solution obtained using quadratic splines with $16 \times 16$ elements is shown in Figure 4a. One can observe that both the pressure and velocity solutions are oscillation-free for all considered cases. Optimal convergence rates are obtained for both the velocity and the pressure field. For the $L^2$-norm and $H^1$-norm of the velocity error, Figure 4b and 4c respectively, asymptotic rates of $k+1$ and $k$ are obtained. For the $L^2$-norm of the pressure shown in Figure 4d we observe asymptotic rates of approximately $k+\frac{1}{2}$, which is half an order higher than those of the $H^1$-norm of the velocity error. For inf-sup compatible discretization pairs where the degrees of the pressure and velocity spaces are $k-1$ and $k$ respectively, the rate of convergence of the $L^2$-norm of the pressure error is known to be equal to that of the $H^1$-norm of the velocity error. We attribute the improved rate for the pressure error using equal order spaces to the fact that compared to the compatible setting the pressure space is one order higher.

---
[2]The bandwidth is defined as the smallest non-negative integer $b$ such that $S_{ij} = 0$ if $|i - j| > b$.

In Figure 5 we study the sensitivity of the computed result with respect to the Skeleton-Penalty stabilization parameter $\gamma$. The $h$-convergence behavior of the solution using $C^1$-continuous quadratic B-splines is studied for a wide range of stabilization parameters, *viz.* $\gamma \in (5 \times 10^{-6}, 1)$. We observe that for this range, the stabilization parameter does not significantly affect the accuracy of the velocity field in the $L^2$-norm and $H^1$-norm, see Figures 5a and 5b, respectively. This is an expected result, as the introduced Skeleton-Penalty term acts only on the pressure field. The pressure solution accuracy is affected by the selection of the stabilization parameter, see Figure 5c. Choosing $\gamma$ too large will lead to a loss of accuracy of the solution, while taking $\gamma$ too small will lead to a loss of stability (this aspect will be discussed in detail below). Figure 5 conveys, however, that the parameter can be selected from a wide range without a significant effect on the accuracy. For the case considered here accuracy deterioration remains very limited in the range $\gamma \in (5 \times 10^{-4}, 5 \times 10^{-2})$. Moreover, for all considered cases we observe the rate of convergence to be independent of the choice of $\gamma$.

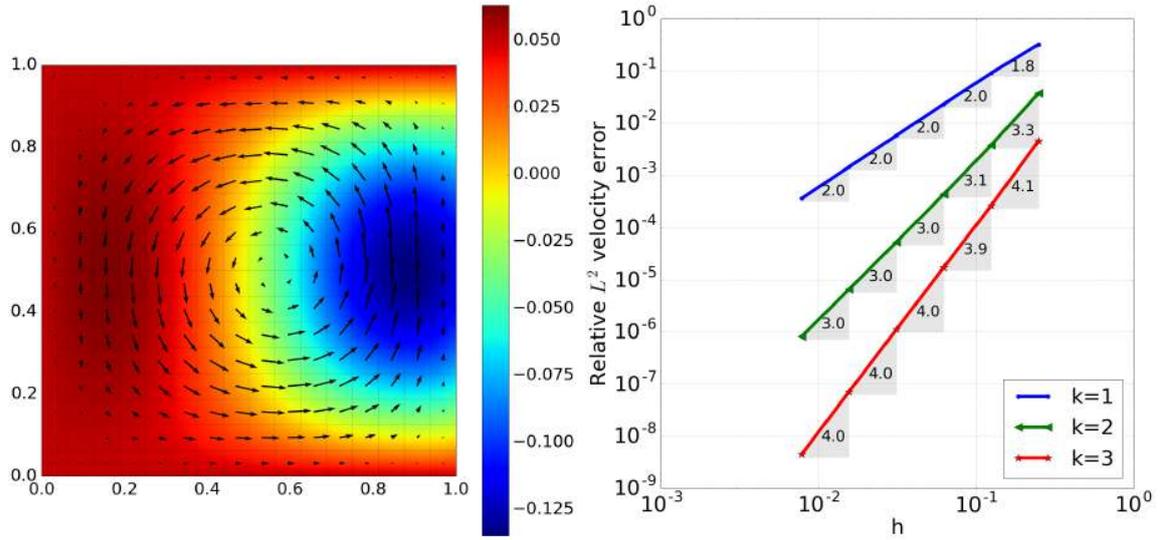

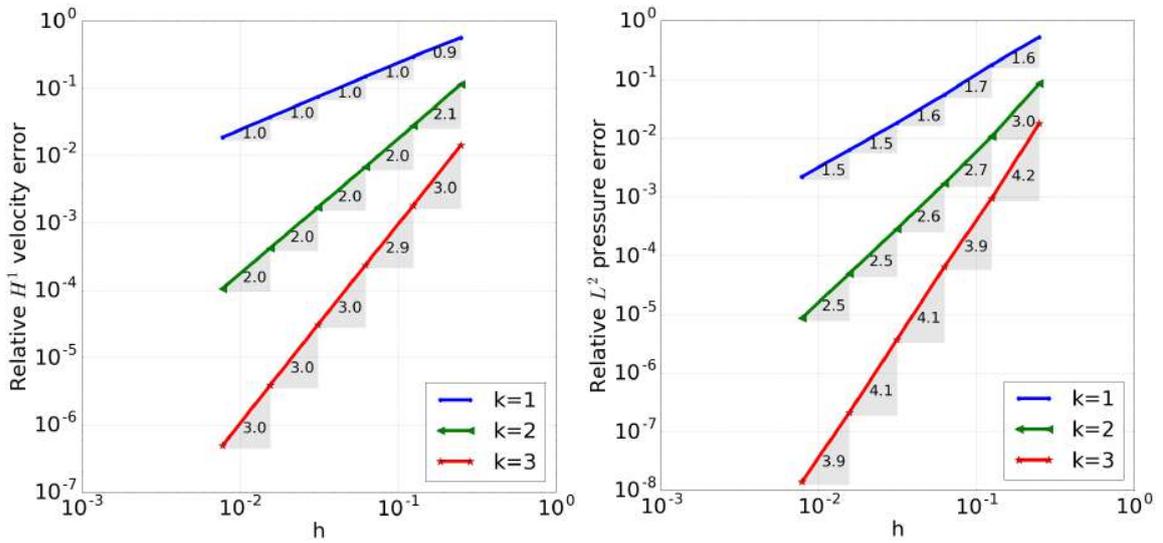

Figure 4: (a) Solution for the steady Stokes problem in Section 5.1, pressure (color) and velocity (vector field). (b-d) Mesh convergence results for B-splines of order $k = 1, 2, 3$ and $C^{k-1}$ regularity.

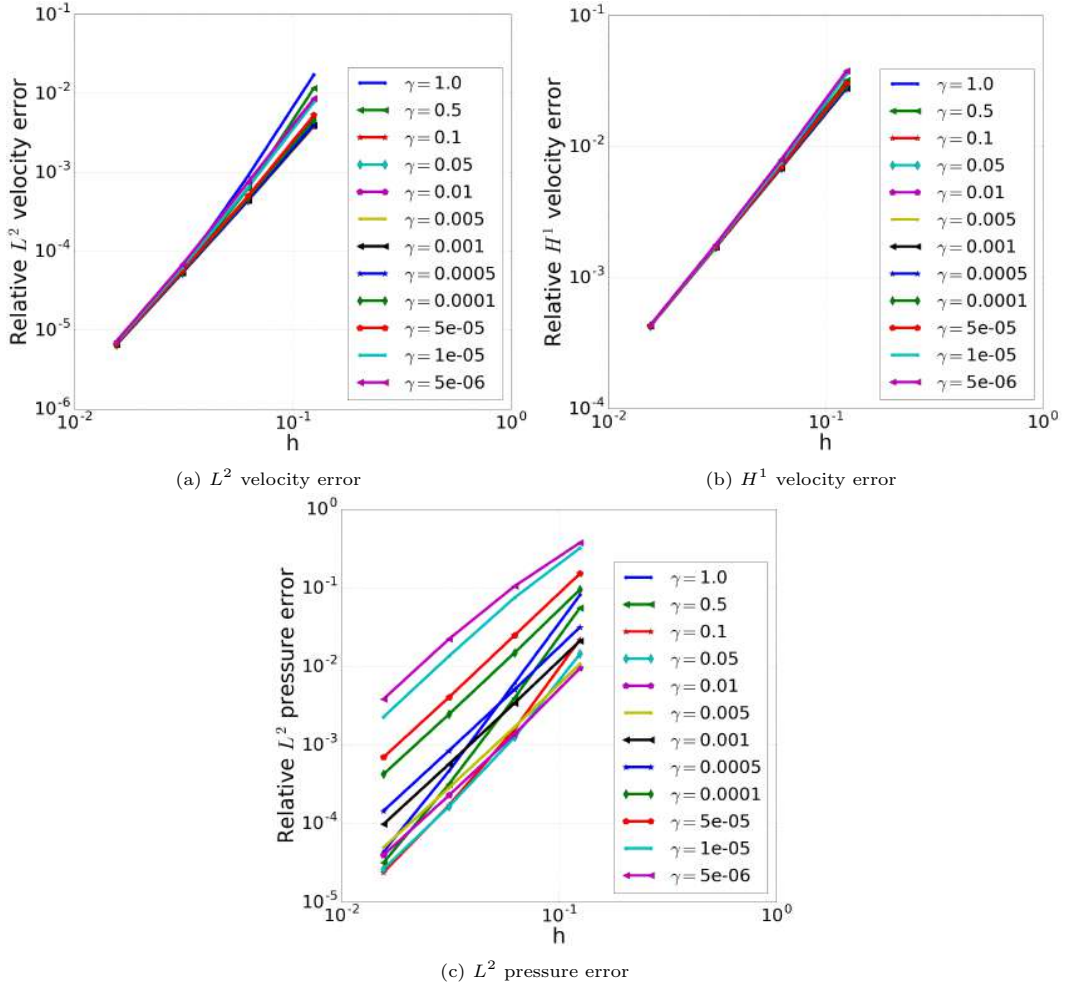

(a) $L^2$ velocity error

(b) $H^1$ velocity error

(c) $L^2$ pressure error

Figure 5: Sensitivity of the quadratic spline approximation of the Stokes problem on the unit square with respect to the stabilization parameter $\gamma$.

To assess the stability of the proposed method, we compute the generalized inf-sup constant (see, e.g., Ref. [23]) associated with the stabilized mixed matrix

$$\begin{bmatrix} \mathbf{A} & \mathbf{B}^T \\ \mathbf{B} & -\mathbf{S} \end{bmatrix},$$

where $\mathbf{A}$, $\mathbf{B}$, and $\mathbf{S}$ are defined as in Sec. 4.2. The discrete stability constant, $\beta_h$, can be computed as the square root of the smallest non-zero eigenvalue of the generalized eigenvalue problem

$$(\mathbf{B}\mathbf{A}^{-1}\mathbf{B}^{\mathbf{T}} + \mathbf{S})\mathbf{q} = \beta_h{}^2 \mathbf{M}_{pp}\mathbf{q}, \tag{35}$$

where $\mathbf{M}_{pp}$ is the Gramian matrix associated with the pressure basis, i.e., $(\mathbf{M}_{pp})_{ij} = (N_i^p, N_j^p)_{Q^h}$. The discrete norm in the pressure space associated with the Gramian matrix is defined as

$$\left\|q^h\right\|_{Q^h}^2 := \left\|q^h\right\|_{L^2(\Omega)}^2 + \gamma \sum_{F \in \mathcal{F}_{skeleton}^h} \int_F \mu^{-1} h_F^{2\alpha+3} \left|[\![\partial_n^{\alpha+1} q^h]\!]\right|^2 \, d\Gamma. \tag{36}$$

Since the norm $\|\cdot\|_{Q^h}$ is stronger than $\|\cdot\|_{L^2(\Omega)}$, numerical inf-sup stability in $\|\cdot\|_{Q^h}$ implies stability for the case that the Gramian matrix $\mathbf{M}_{pp}$ is defined as the $L^2$ pressure mass matrix.

Figure 6 present the results of the numerical stability study conducted for various selections of the stabilization parameter $\gamma$. The results convey that the proposed method is stable evidenced by the fact that the discrete stability constants are bounded from below away from zero under mesh refinement. Figure 6a presents the results for the order-dependent choice of $\gamma$ considered above. For this choice the stability parameter is virtually independent of the order of the approximation, and optimal rates of convergence for the $L^2$ pressure error are obtained. Choosing $\gamma$ too large ( Figure 6b) does not affect the stability of the formulation, but negatively affects the accuracy of the solutions. Choosing $\gamma$ too small ( Figure 6c) does affect the stability in the sense that the discrete generalized inf-sup constant reduces with increasing degree $k$, although it is still essentially independent of $h$. The reduced stability evidently also affects the accuracy of the solution. The results for this choice of $\gamma$ reflect that the higher-order regularity of B-splines has a positive effect on the stability, in the sense that for the same selection of the stability parameter, the discrete generalized inf-sup constant increases with increasing regularity.

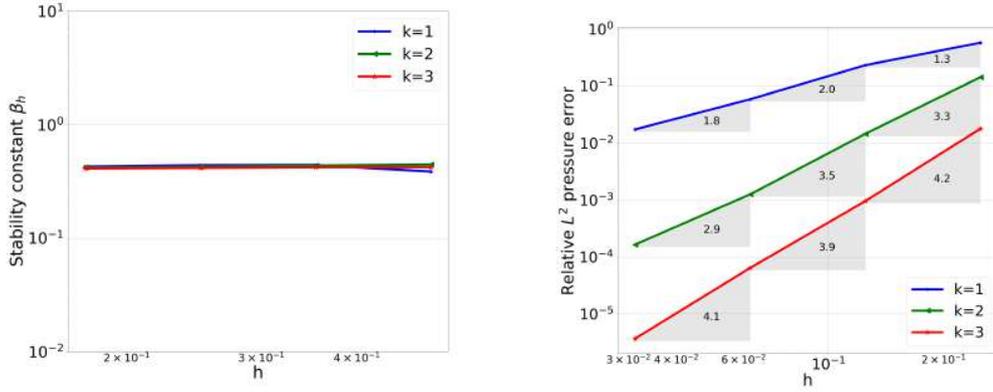

(a) $\gamma = 1$ ($k = 1$), $5 \times 10^{-2}$ ($k = 2$), $10^{-3}$ ($k = 3$)

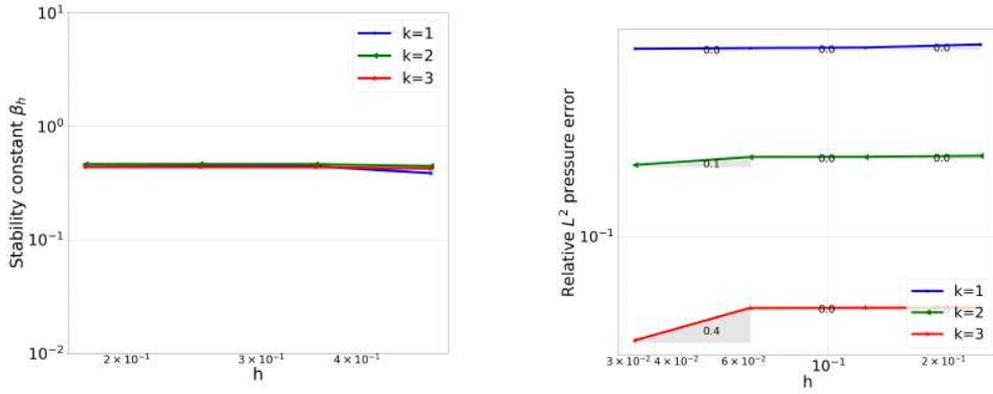

(b) $\gamma = 10^5$ (too large)

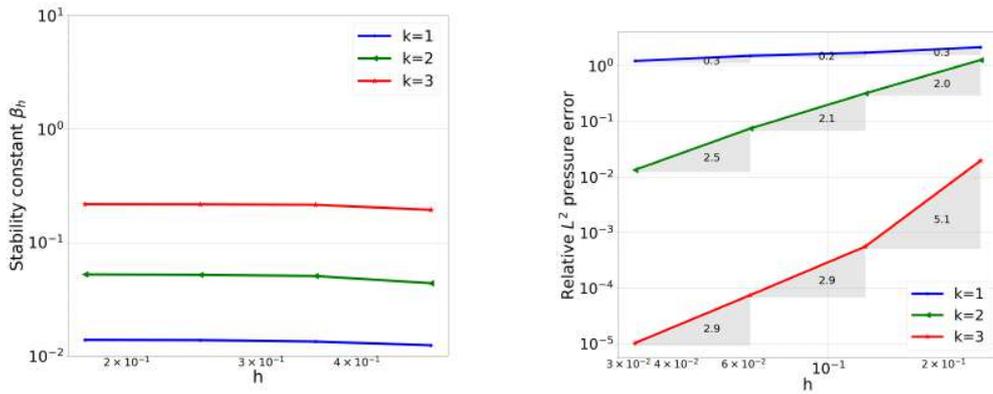

(c) $\gamma = 10^{-5}$ (too small)

Figure 6: Numerical study of the discrete stability constant (left column) and the $L^2$ pressure error (right column) for various mesh sizes ($h$) and spline degrees ($k = 1, 2, 3$) of highest regularities. Figure 6a corresponds to the stabilization parameter choice $\gamma = 1$ ($k = 1$), $5 \times 10^{-2}$ ($k = 2$), $10^{-3}$ ($k = 3$), as considered above. Figure 6b considers the case for which the stabilization parameter is chosen too large, whereas Figure 6c represents a too small selection of the stabilization parameter.

The performance of the proposed Skeleton-stabilized IsoGeometric Analysis framework is further studied based on the generalized Stokes equations with homogeneous Dirichlet boundary conditions:

$$\begin{cases} \text{Find } \mathbf{u}: \overline{\Omega} \to \mathbb{R}^d, \text{ and } p: \overline{\Omega} \to \mathbb{R} \text{ such that:} \\ \sigma\mathbf{u} - \nabla \cdot (2\mu\nabla^s\mathbf{u}) + \nabla p = \mathbf{f} \quad \text{in } \Omega, \\ \qquad\qquad\qquad\qquad \nabla \cdot \mathbf{u} = 0 \quad \text{in } \Omega, \\ \qquad\qquad\qquad\qquad\qquad \mathbf{u} = \mathbf{0} \quad \text{on } \Gamma_D, \end{cases} \qquad (37)$$

This system – for which the body force $\mathbf{f}$ is selected in accordance with the manufactured solution (34) – is characterized by the Damköhler number

$$\text{Da} = \frac{\sigma L^2}{\mu}, \qquad (38)$$

where $\sigma$ is the reaction coefficient, and $L$ is a characteristic length scale for the problem (in this case the width/height of the unit square).

In Figure 7 we study the $h$-convergence behavior of $C^{k-1}$-continuous B-splines for various degrees $k = 1, 2, 3$ and $\text{Da} = 1, 10, 1000$. To control the reaction term, we supplement the stabilization term with a contribution from $\sigma$ to the scaling ratio, *i.e.*,

$$s(p^h, q^h) = \sum_{F \in \mathcal{F}^h_{skeleton}} \int_F \gamma(\mu + \sigma h_F^2)^{-1} h_F^{2\alpha+3} [\![\partial_n^{\alpha+1} p^h]\!] [\![\partial_n^{\alpha+1} q^h]\!] d\Gamma.$$

The stabilization parameter is now chosen equal to $\gamma = 1$ ($k = 1$), $5e - 2$ ($k = 2$), $1e - 3$ ($k = 3$). Note that the non-reactive case of $\text{Da} = 0$ corresponding to $\sigma = 0$ resembles the case considered above. For all considered cases we observe the approximation of the velocity solution and pressure solution to be virtually independent of the Damköhler number.

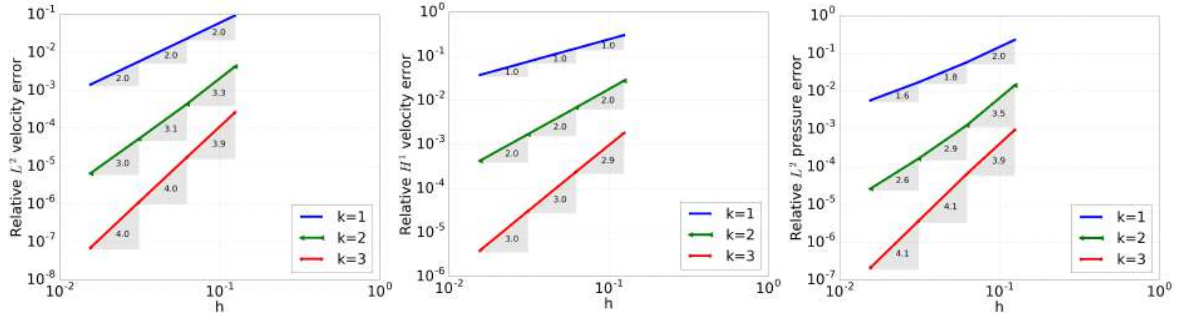

(a) Da = 1

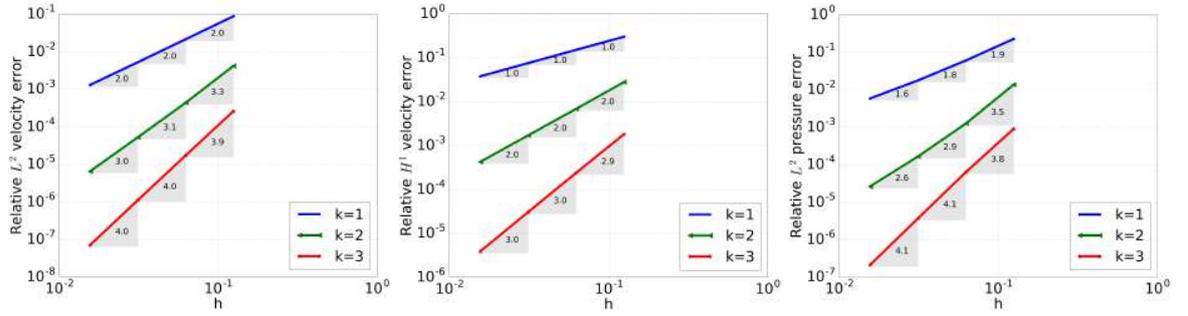

(b) Da = 10

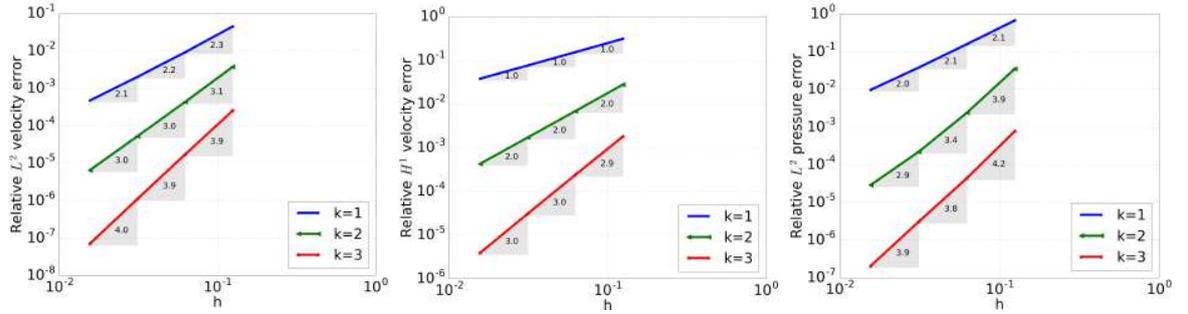

(c) Da = 1000

Figure 7: $h$-convergence behavior of $C^{k-1}$-continuous B-spline spaces of degree $k = 1, 2, 3$ for various Damköhler numbers.

To understand the effect of reduced regularity – which is of particular importance in the case of multi-patch models – we first study the B-spline discretization of the Stokes problem on the unit square with varying intra-patch regularities. That is, we consider the spline discretizations $S_\alpha^k$ of order $k$ with regularity $\alpha = 0, \ldots, k-1$. A stabilization parameter of $\gamma = 10^{-\alpha}k^{-4}$ – which effectively decreases the penalty parameter with increasing order and regularity – was found to yield an adequate balance between accuracy and stability for the considered simulations. Derivation of a rigorous selection criterion for the penalty parameter is beyond the scope of the current work. Note that because the case of $k = 1$ and $\alpha = 0$ has already been considered above, we here restrict ourselves to the spline degrees $k = 2, 3, 4$. The $h$-convergence results are collected in Figure 8. Note that we plot the errors versus the square root of the number of degrees of freedom to enable comparison of the various approximations. We observe optimal convergence rates for both the velocity and the pressure approximation for all cases. As anticipated the accuracy per degree of freedom improves with increasing regularity. Note that in the case of $S_0^k$ – which is equivalent to the Lagrange basis – we observe similar approximation behavior as for the continuous interior-penalty method [17].

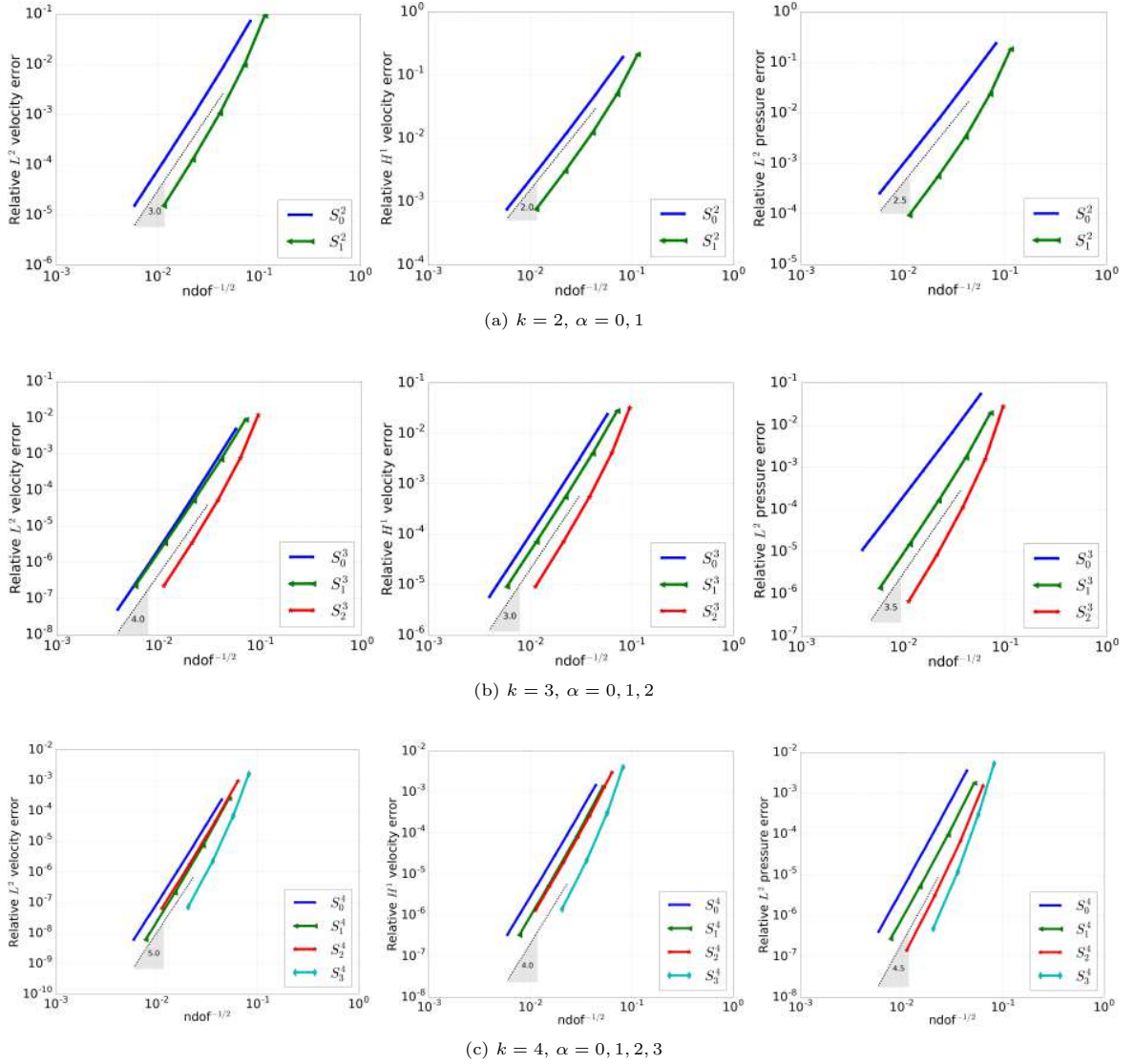

Figure 8: $h$-convergence results for the Stokes problem in a unit square using B-splines spaces $S_\alpha^k$ of various degrees $k = 2, 3, 4$ and regularities $0 \leq \alpha \leq k - 1$.

The stability study based on the discrete generalized inf-sup constant for cases with reduced regularities is shown in Figure 9. For all considered cases the discrete stability constants are observed to be bounded from below away from zero under mesh refinement. Figure 9 presents the results for the case of quadratic, cubic, and quartic splines ($k = 2, 3, 4$), with various orders of regularity ($0 \leq \alpha \leq k - 1$), using the above-mentioned $(k, \alpha)$-dependent stabilization parameter, the inf-sup constant is observed to be virtually independent of the regularity. As was also observed in Figure 9, increased regularity enhances the stability in the sense that a smaller stabilization parameter suffices.

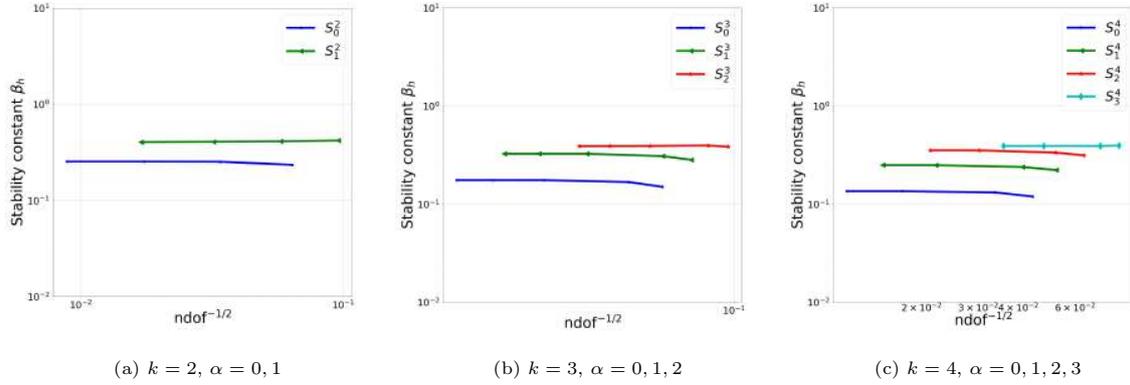

Figure 9: Discrete stability constant and its behavior under mesh refinement for the proposed method using B-splines spaces $S_\alpha^k$ of various degrees $k = 2, 3, 4$ and regularities $0 \leq \alpha \leq k - 1$.

## 5.2. Steady Stokes flow in a quarter annulus ring

To demonstrate the performance of the proposed Skeleton-Penalty stabilization in the context of IsoGeometric Analysis, we consider the steady Stokes problem in the open quarter annulus domain

$$\Omega = \left\{ \mathbf{x} \in \mathbb{R}^2_{>0} : R_1 < |\mathbf{x}| < R_2 \right\},$$

with inner radius $R_1 = 1$ and outer radius $R_2 = 4$. We parametrize this domain using NURBS. Homogeneous Dirichlet boundary conditions are prescribed on the entire boundary $\partial\Omega = \Gamma_D$ and, accordingly, it holds that $\Gamma_N = \emptyset$. The body force $\mathbf{f}$ is selected in accordance with the manufactured solution [24, 14]

$$\mathbf{u}(\mathbf{x}) = \begin{pmatrix} 10^{-6} x^2 y^4 (x^2 + y^2 - 1)(x^2 + y^2 - 16)(5x^4 + 18x^2 y^2 - 85x^2 + 13y^4 - 153y^2 + 80) \\ 10^{-6} xy^5 (x^2 + y^2 - 1)(x^2 + y^2 - 16)(102x^2 + 34y^2 - 10x^4 - 12x^2 y^2 - 2y^4 - 32) \end{pmatrix}, \quad (39a)$$

$$p(\mathbf{x}) = 10^{-7} xy(y^2 - x^2)(x^2 + y^2 - 16)^2 (x^2 + y^2 - 1)^2 \exp\left(14(x^2 + y^2)^{-1/2}\right). \quad (39b)$$

Note that $\mathbf{u}$ vanishes on $\partial\Omega$ in accordance with the Dirichlet boundary condition. Moreover, the pressure complies with $\int_\Omega p \, d\Omega = 0$. This manufactured solution is illustrated in Figure 10.

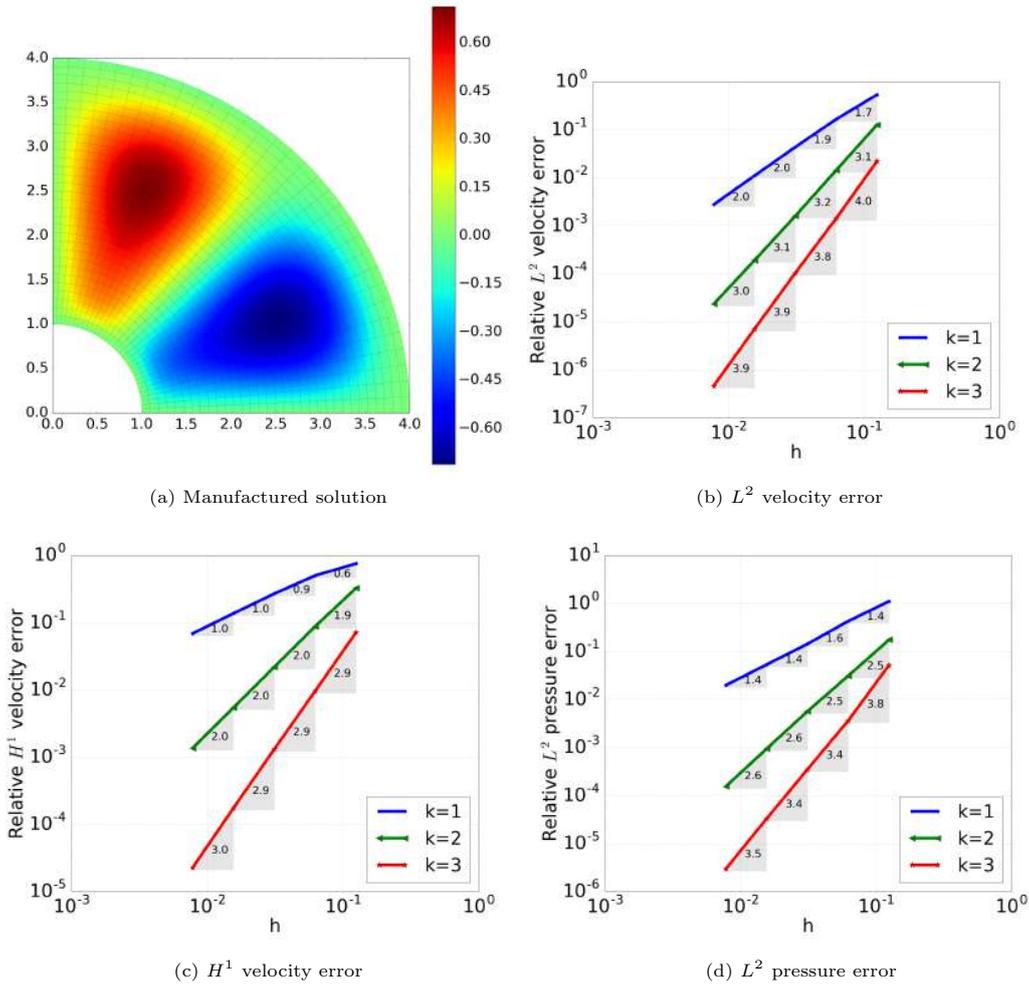

Figure 10: (a) Pressure solution for the steady Stokes problem on a quarter annulus ring in Section 5.2. (b-d) $h$-convergence results for B-splines of order $k = 1, 2, 3$ and $C^{k-1}$ regularity.

In this example we consider B-spline bases of orders $k = 1, 2, 3$ on meshes ranging from $8 \times 8$ to $128 \times 128$ elements. We divert here from the isoparametric concept in order to also study the performance of linear bases, which are incapable of parametrizing the annulus ring exactly. We will consider NURBS-based isogeometric analysis in later test cases. For the simulation, the stabilization parameter is taken as $\gamma = 1$ ($k = 1$), $5 \cdot 10^{-2}$ ($k = 2$), $1 \cdot 10^{-3}$ ($k = 3$). In Figure 10a the pressure solution obtained using $C^1$-continuous quadratic B-splines on a $32 \times 32$ element mesh is shown, which is observed to be free of oscillations. In Figures 10b and 10c we observe optimal convergence rates for the velocity error of $k+1$ for the $L^2$-norm and $k$ for the $H^1$-norm, respectively. As for the unit square problem considered above, an asymptotic rate of convergence of approximately $k + \frac{1}{2}$ is observed for the $L^2$-norm of the pressures.

### 5.3. Steady Navier-Stokes flow in a full annulus domain

As the baseline test case for the Skeleton-stabilized IsoGeometric analysis of the steady incompressible Navier-Stokes equations we consider the cylindrical Couette flow between two cylinders as shown in Figure 11a, which was studied in the context of compatible spline discretizations in [10]. The outer cylinder is fixed, while the inner cylinder rotates with surface velocity $U = \omega R_1$. For low Reynolds numbers the flow in between the cylinders will remain steady, two-dimensional, and axisymmetric. The analytical velocity solution of the problem is then given by

$$u = \begin{pmatrix} -(Ar + Br^{-1})\sin(\theta) \\ (Ar + Br^{-1})\cos(\theta) \end{pmatrix}, \tag{40}$$

where $(r, \theta)$ are the polar coordinates originating from the center of the cylinders, and

$$A = -U\frac{\delta^2}{R_1(1-\delta^2)}, \qquad B = U\frac{R_1}{1-\delta^2}, \tag{41}$$

with $\delta = R_1/R_2$ the ratio of radii of the inner and outer cylinders. The analytical pressure solution is a constant function, which supplemented with the zero average pressure condition $\int_\Omega p\, d\Omega = 0$ results in a zero pressure field. Here we consider the case of $\omega = 1$, $R_1 = 1$, and $R_2 = 2$. The solution for this case is illustrated in Figure 11c.

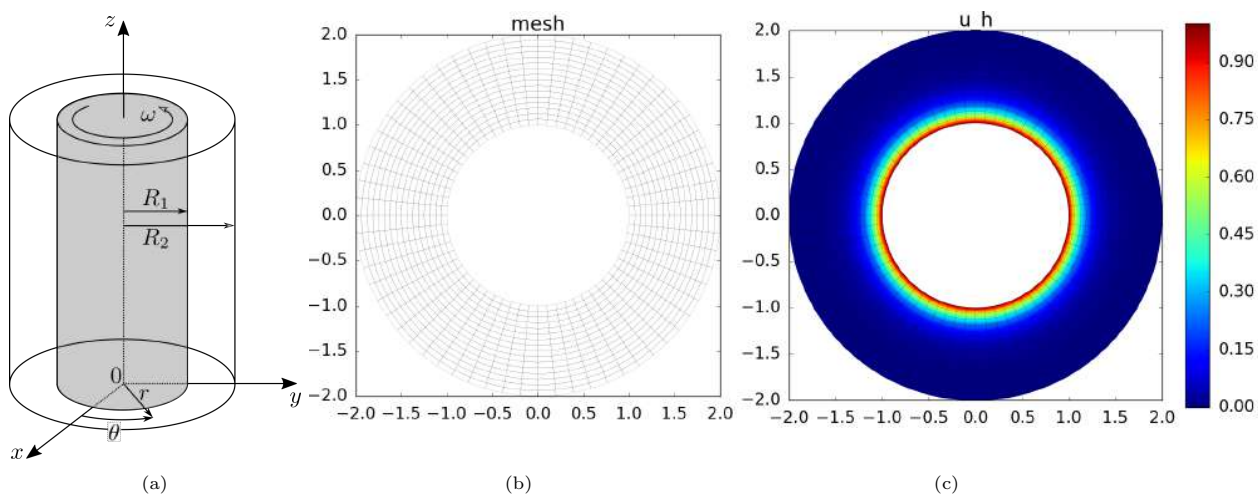

Figure 11: (a) Setup of the cylindrical Couette flow problem. (b) Two-dimensional polar mesh, and (c) a typical solution of the radial velocity component.

For the parametrization of the geometry the polar map

$$(0,1)^2 \ni (\xi_1, \xi_2) \mapsto \mathbf{F}(\xi_1, \xi_2) = \begin{pmatrix} ((R_2 - R_1)\xi_2 + R_1)\sin(2\pi\xi_1) \\ ((R_2 - R_1)\xi_2 + R_1)\cos(2\pi\xi_1) \end{pmatrix} \tag{42}$$

is used, where $(\xi_1, \xi_2)$ are the coordinates of the unit square parameter domain. The problem is discretized using B-splines of degree $k = 1, 2, 3$ with $C^{k-1}$-regularity, which are periodic in the circumferential $\xi_1$-direction. In Figure 12 we study the mesh convergence behavior of the velocity approximation in the $L^2$-norm and $H^1$-norm. The coarsest mesh considered consists of $8 \times 2$ elements (two elements in the radial direction), which is uniformly refined until a mesh of $128 \times 32$ elements is obtained. We observed optimal rates of convergence for all orders in both the $L^2$-norm and $H^1$-norm. The pre-asymptotic behavior observed for the $H^1$-norm is a result of the fact that the boundary layer near the inner circle is not even remotely resolved by a single element. By virtue of the nature of the problem, the analytical zero pressure field is satisfied identically.

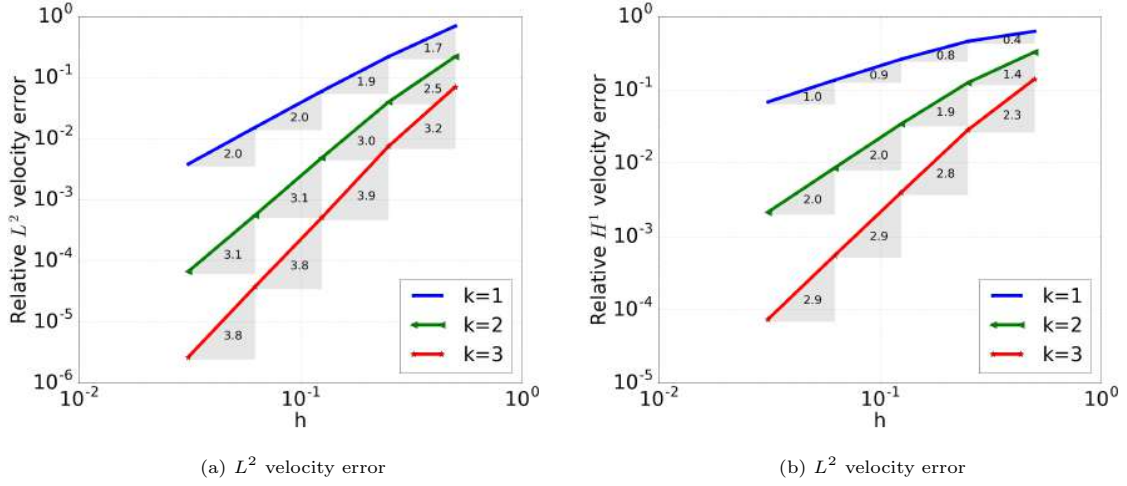

(a) $L^2$ velocity error

(b) $L^2$ velocity error

Figure 12: $h$-convergence study of the cylindrical Couette flow problem using various order B-splines with $C^{k-1}$ regularity.

5.4. Navier-Stokes flow around a circular cylinder

To study the performance of the proposed formulation for the Navier-Stokes equations in further detail we consider the benchmark problem proposed by Schäfer and Turek [25]. In this benchmark the flow around a cylinder which is placed in a channel is studied. The geometry of this test case is shown in Figure 13, where the channel length is $L = 2.2$ m, the channel height is $H = 0.41$ m, and the cylinder radius is $R = 0.05$ m. The center of the cylinder is positioned at $\frac{1}{2}(W, W) = (0.2, 0.2)$ m, which has an offset of $\frac{1}{2}\delta = 0.005$ m with respect to the center line of the channel (such that $W = H - \delta = 0.4$ m). At the inflow boundary ($x = 0$) a parabolic horizontal flow profile is imposed

$$\mathbf{u}(0, y) = \begin{pmatrix} 4U_m y(H - y)/H^2 \\ 0 \end{pmatrix}$$

with maximum velocity $U_m$. A no slip boundary condition is imposed along the bottom and top boundaries, as well as along the surface of the cylinder. At the outflow boundary ($x = L$) a zero traction boundary condition is used. The density and kinematic viscosity of the fluid are taken as $\rho = 1.0 \, \text{kg/m}^3$ and $\mu = 1 \times 10^{-3} \, \text{m}^2/\text{s}$, respectively.

We consider two cases, one corresponding to an inflow velocity that results in a steady flow, and one corresponding to an inflow velocity that results in an unsteady flow. These two cases are characterized by the Reynolds number

$$\text{Re} = \frac{2\bar{U}R}{\mu},$$

where $\bar{U} = \frac{2}{3}U_m$ is the mean inflow velocity. As quantities of interest we consider the drag and lift coefficients

$$c_D = \frac{F_D}{\rho \bar{U}^2 R}, \qquad c_L = \frac{F_L}{\rho \bar{U}^2 R},$$

where $F_D$ and $F_L$ are the resultant lift and drag forces acting on the cylinder. These forces are weakly evaluated as (see e.g., [26])

$$F_D = \mathcal{R}(\mathbf{u}, p; \boldsymbol{\ell}_1), \qquad F_L = \mathcal{R}(\mathbf{u}, p; \boldsymbol{\ell}_2),$$

where

$$\mathcal{R}(\mathbf{u}, p; \boldsymbol{\ell}_i) := (\partial_t \mathbf{u}, \boldsymbol{\ell}_i) + c(\mathbf{u}; \mathbf{u}, \boldsymbol{\ell}_i) + a(\mathbf{u}, \boldsymbol{\ell}_i) + b(p, \boldsymbol{\ell}_i),$$

with $\boldsymbol{\ell}_i \in [H^1_{0,\partial\Omega\setminus\Gamma}(\Omega)]^d$ and $\boldsymbol{\ell}_i|_\Gamma = -\mathbf{e}_i$, $i = 1, 2$. We note that these lift and drag evaluations are consistent with the weak formulation (18), and are different from the formulations given in [27] and [28] where in the former, the time derivative term is neglected (so only consistent for the steady case), and in the latter, both the convective term and time derivative term are neglected (thus only consistent for the case of steady Stokes equations). For the steady test case, we also consider the pressure drop over the cylinder

$$\Delta p = p(W/2 - R, W/2) - p(W/2 + R, W/2),$$

and for the unsteady test case, we consider the Strouhal number

$$St = \frac{Df}{\bar{U}}$$

as additional quantities of interest, where $f$ is the frequency of vortex shedding and $D$ is the diameter of the cylinder.

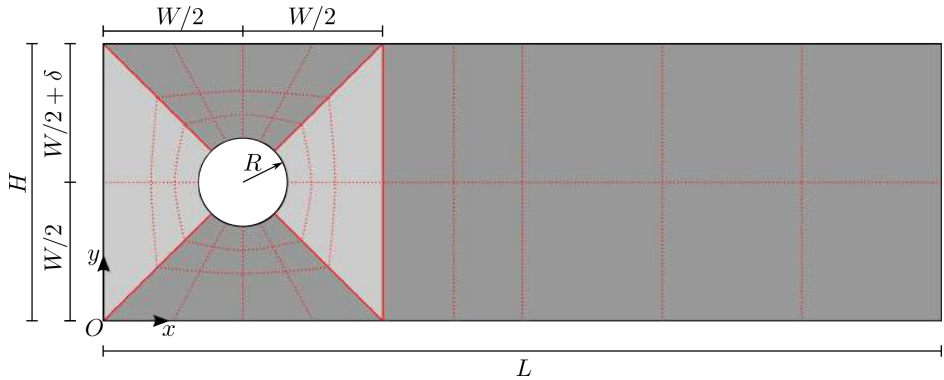

Figure 13: Multi-patch parametrization of the channel flow problem with a circular obstacle.

The geometry is parameterized by a quadratic ($k = 2$) multi-patch NURBS surface, as shown schematically in Figure 13. The boundaries between the five patches are indicated by the solid red lines, while the element boundaries within the patches are marked by dashed red lines. Full $C^{k-1}$-continuity is maintained at the intra-patch element boundaries. For the coarsest mesh we employ $8 \times 5$ elements in the circumferential and radial direction, respectively, for each of the four patches adjacent to the cylinder. The discretization of the downstream patch conforms with its neighboring patch and consists of $8 \times 8$ elements in the vertical and horizontal direction, respectively. The employed NURBS are non-uniform as the meshes are locally refined toward the cylinder, and coarsened toward the outflow boundary.

We discretize both the velocity components and the pressure using the NURBS basis employed for the geometry parametrization, making this a true isogeometric analysis. Our coarsest quadratic NURBS mesh is refined uniformly to study the $h$-convergence behavior of the above-mentioned quantities of interest. Moreover, we elevate the order of our coarsest mesh to a cubic ($k = 3$) multi-patch NURBS surface with $C^{k-1}$-continuity inside the patches, and subsequently perform uniform mesh refinements to study the $h$-convergence behavior for the cubic case.

*5.4.1. Steady flow*

We first consider the case of Reynolds number Re = 20, for which a steady flow develops. The velocity magnitude and pressure solutions for this case are shown in Figure 14. A two times uniform refinement of the coarsest quadratic NURBS mesh is used to compute this result, which contains $n_{\text{dof}} = 12180$ degrees of freedom. The computed drag and lift coefficients, $c_D = 5.5798$ and $c_L = 0.010605$, are in excellent agreement, respectively, with the benchmark ranges $(5.57, 5.59)$, and $(0.0104, 0.0110)$ reported in [25], as is the computed solution for the pressure drop $\Delta p = 0.117514$.

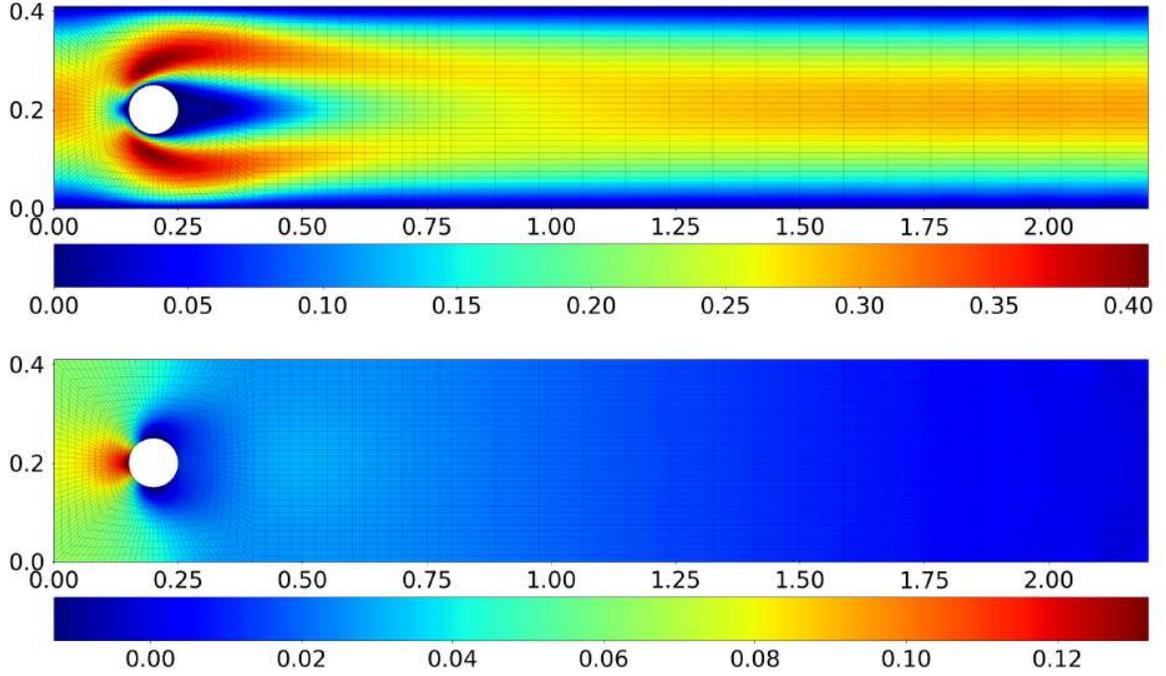

Figure 14: Velocity magnitude (top) and pressure (bottom) solutions of the steady cylinder flow problem using quadratic NURBS with $n_{\text{dof}} = 12180$.

In Figure 15 we present the $h$-convergence results for the three quantities of interest. For the quadratic case in Figure 15a we consider five meshes, where the coarsest one corresponding to the geometry parameterization, results in 1056 degrees of freedom, and the four times uniformly refined mesh results in 177636 degrees of freedom. The errors are computed with respect to the high-quality reference values proposed in [29]:

$$C_D^{ref} = 5.57953523384, \qquad C_L^{ref} = 0.010618948146, \qquad \Delta_p^{ref} = 0.11752016697$$

We observe convergence of all three quantities of interest to the benchmark solutions. In particular for the lift coefficient and the pressure drop the observed asymptotic rates match well with the expected optimal rates of $2k$ [30]. In Figure 15b we consider the mesh convergence of the quantities of interest for the cubic NURBS case, for which the coarsest mesh consists of 1356 degrees of freedom, and the finest mesh (4 uniform refinements) consists of 181356 degrees of freedom. As expected we observed improved rates of convergence compared to the quadratic case. Note that in terms of degrees of freedom there is virtually no difference between the finest quadratic mesh and the finest cubic mesh, which conveys that increasing the spline order is favorable from an accuracy per degree of freedom point of view. We expect that the irregular behavior of the convergence rate for the lift coefficient on the finest cubic meshes is related to approaching the accuracy of the reference solution from [29].

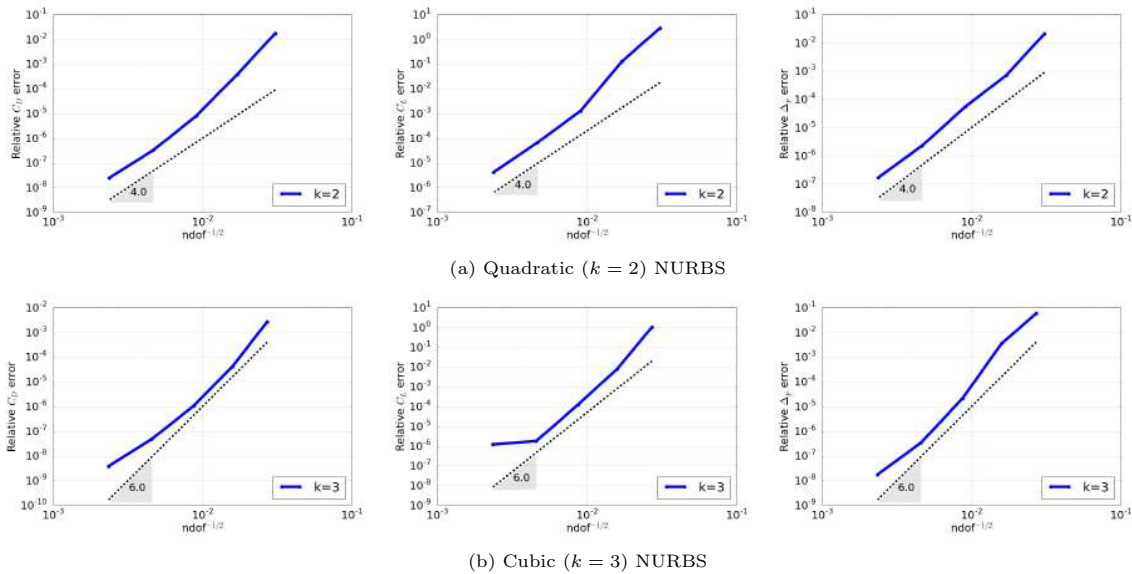

Figure 15: $h$-convergence results for the drag coefficient (left column), lift coefficient (middle column) and pressure drop (right column) of the steady cylinder flow problem for quadratic (top row) and cubic (bottom row) NURBS.

#### 5.4.2. Unsteady flow

For the case of Reynolds number Re = 100 there is no longer a steady solution. Instead, once the flow is fully developed, oscillatory vortex shedding occurs, as illustrated by the snapshot shown in Figure 16. For this figure, a two times uniformly refined quadratic NURBS parametrization is used, which results in a total of 12180 degrees of freedom. In order to capture the vortex shedding, the downstream mesh characteristics have been adjusted in comparison to the steady test case, in the sense that the refinement zone stretches out further behind the cylinder. We have used a time step of $\Delta t = 1/20$ s for the first 4 s of the simulations in order to let the flow develop, after which we switch to a smaller time step size of $\Delta t = 1/200$ s to accurately capture the oscillatory behavior of the solution. In Figure 17 the evolution of the drag coefficient, lift coefficient and pressure drop over time is shown for the fully developed vortex shedding flow.

Table 1 presents a comparision result for three consecutive uniform mesh refinement levels using quadratic NURBS and $\Delta t = 1/200$ s. The flow is only considered when it is fully developed. The time cycle is arbitrarily chosen such that at the start and end of the interval, the lift coefficients attain two consecutive local minima. The quantities of interest are the minimum and maximum of the lift and drag coefficients, the length of the time cycle, and the Strouhal number. From Table 1, we compute Table 2, which shows the relative errors of the quantities of interest (and their convergence rates). We observe that these quantities of interest converge very well to the high-quality results reported in [31]. At the first level of refinement with only 3420 degrees of freedom, the results already start to be close to the reference values, with the relative errors of $5.53 \times 10^{-3}$ and $8.96 \times 10^{-3}$ for the minimum and maximum of the drag coefficient, and approximately $8 \times 10^{-2}$ for the minimum and maximum of the lift coefficient. At the third level of refinement with 45828 degrees of freedom, when the mesh is fine enough to resolve the boundary layer around the cylinder, and to accurately capture the dynamics of the flow, we obtain the convergence rates of $2k$ ($k = 2$) as in the steady test case, with errors of $1.34 \times 10^{-4}$ and $8.08 \times 10^{-5}$ for the minimum and maximum of the drag coefficient, and $1.49 \times 10^{-3}$ and $1.12 \times 10^{-3}$ for the minimum and maximum of the lift coefficient, respectively. The obtained time cycles and Strouhal numbers are also in a good agreement with the reference values. Note that the computation of these quantities is based on the time interval of two consecutive local minima of the lift coefficient, and are therefore directly affected by the choice of time step.

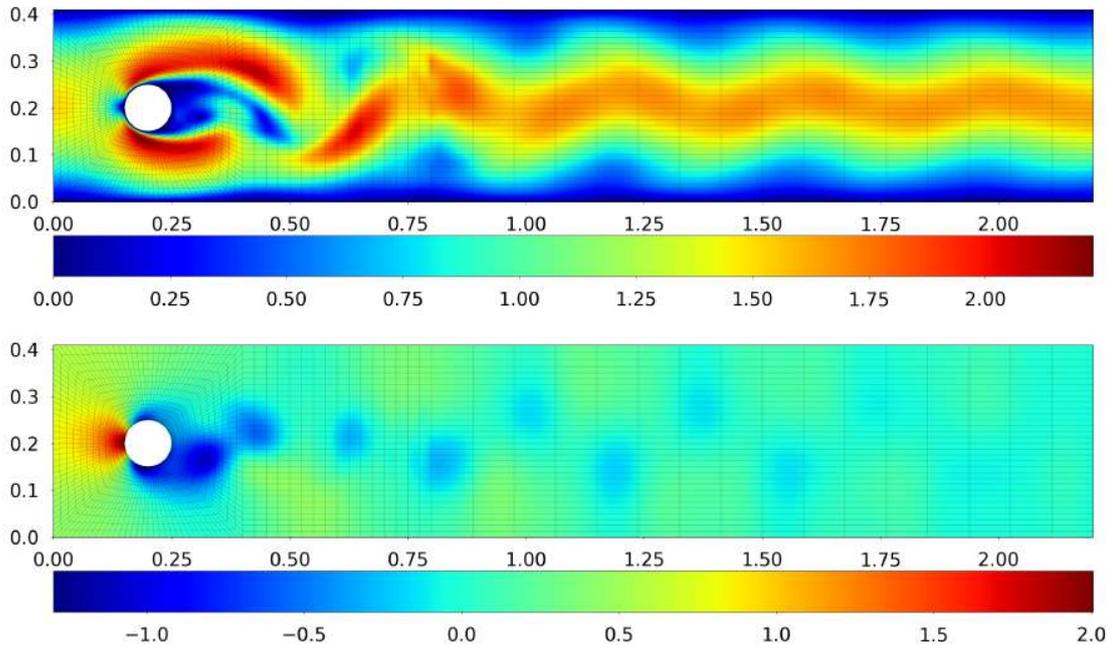

Figure 16: A snapshot of velocity (top) and pressure (bottom) of the unsteady cylinder flow problem (Re=100); A Von Kármán vortex street is clearly visible behind the cylinder.

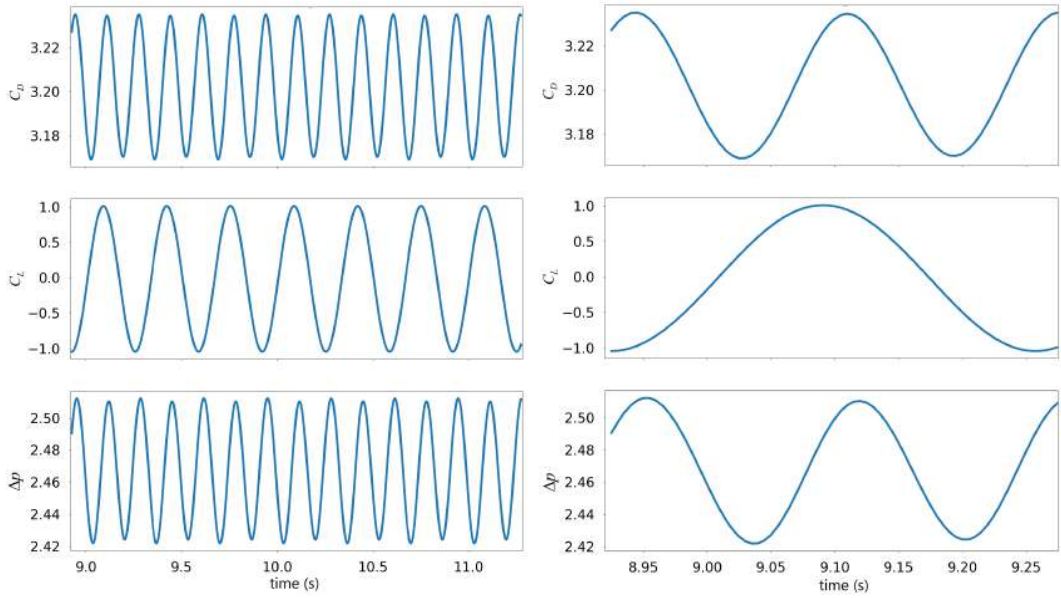

Figure 17: Drag coefficient, lift coefficient and pressure drop over time (left) and a zoom of one period (right) for the unsteady cylinder flow problem. These results are based on a quadratic NURBS $k = 2$ discretization with two levels of refinements from the coarsest mesh.

Table 1: Minimum and maximum of the drag and lift coefficients, time cycle length, and the Strouhal number for the unsteady cylinder flow problem. For all cases the degree is $k = 2$ and the time step size is $\Delta t = 1/200$.

| Level | $n_{\text{dof}}$ | min $C_D$ | max $C_D$ | min $C_L$ | max $C_L$ | $1/f$ | St |
|---|---|---|---|---|---|---|---|
| 1 | 3420 | 3.18175 | 3.25632 | -1.10535 | 1.06472 | 0.34500 | 0.28986 |
| 2 | 12180 | 3.16893 | 3.23507 | -1.04482 | 1.00895 | 0.33500 | 0.29851 |
| 3 | 45828 | 3.16469 | 3.22765 | -1.01977 | 0.98547 | 0.33000 | 0.30303 |
| Ref [31]: 6 | 667264 | 3.16426 | 3.22739 | -1.02129 | 0.98657 | 0.33125 | 0.30189 |

Table 2: Relative error convergence of the minimum and maximum of the drag and lift coefficients of the unsteady cylinder flow problem, computed from Table 1. For all cases the degree is $k = 2$ and the time step size is $\Delta t = 1/200$. The rate of convergence is here indicated by $r$.

| Level | $n_{\text{dof}}$ | error min $C_D$ | error max $C_D$ | error min $C_L$ | error max $C_L$ |
|---|---|---|---|---|---|
| 1 | 3420 | $5.53 \times 10^{-3}$ | $8.96 \times 10^{-3}$ | $8.23 \times 10^{-2}$ | $7.92 \times 10^{-2}$ |
| 2 | 12180 | $1.47 \times 10^{-3}$ ($r = 2.08$) | $2.38 \times 10^{-3}$ ($r = 2.09$) | $2.30 \times 10^{-2}$ ($r = 2.00$) | $2.26 \times 10^{-2}$ ($r = 1.97$) |
| 3 | 45828 | $1.34 \times 10^{-4}$ ($r = 3.62$) | $8.08 \times 10^{-5}$ ($r = 5.11$) | $1.49 \times 10^{-3}$ ($r = 4.14$) | $1.12 \times 10^{-3}$ ($r = 4.54$) |

### 5.5. Three-dimensional Navier-Stokes flow in a sphere

To demonstrate the performance of the Skeleton-Penalty formulation in the three-dimensional case, we consider the 3D benchmark problem of Navier-Stokes flow proposed by Ethier and Steinman in [32] with the domain considered a sphere. We parametrize the spherical geometry by mapping a bi-unit cube parameter domain $\hat{\Omega} = (-1, 1)^3 \ni \boldsymbol{\xi}$ onto the physical domain $\Omega \ni \mathbf{x}$ through

$$\mathbf{x} = \begin{pmatrix} \xi_1 \sqrt{1 - \frac{\xi_2^2}{2} - \frac{\xi_3^2}{2} + \frac{\xi_2^2 \xi_3^2}{3}} \\ \xi_2 \sqrt{1 - \frac{\xi_3^2}{2} - \frac{\xi_1^2}{2} + \frac{\xi_3^2 \xi_1^2}{3}} \\ \xi_3 \sqrt{1 - \frac{\xi_1^2}{2} - \frac{\xi_2^2}{2} + \frac{\xi_1^2 \xi_2^2}{3}} \end{pmatrix}. \tag{43}$$

We consider the manufactured solution

$$\mathbf{u}(\mathbf{x}) = \begin{pmatrix} -a[e^{ax} \sin(ay + dz) + e^{az} \cos(ax + dy)] \\ -a[e^{ay} \sin(az + dx) + e^{ax} \cos(ay + dz)] \\ -a[e^{az} \sin(ax + dy) + e^{ay} \cos(az + dx)] \end{pmatrix}, \tag{44a}$$

$$p(\mathbf{x}) = -\frac{a^2}{2} \big[ e^{2ax} + e^{2ay} + e^{2az} + 2\sin(ax + dy)\cos(az + dx)e^{a(y+z)} \\ + 2\sin(ay + dz)\cos(ax + dy)e^{a(z+x)} + 2\sin(az + dx)\cos(ay + dz)e^{a(x+y)} \big]. \tag{44b}$$

with parameters $a = 1$ and $d = 1$.

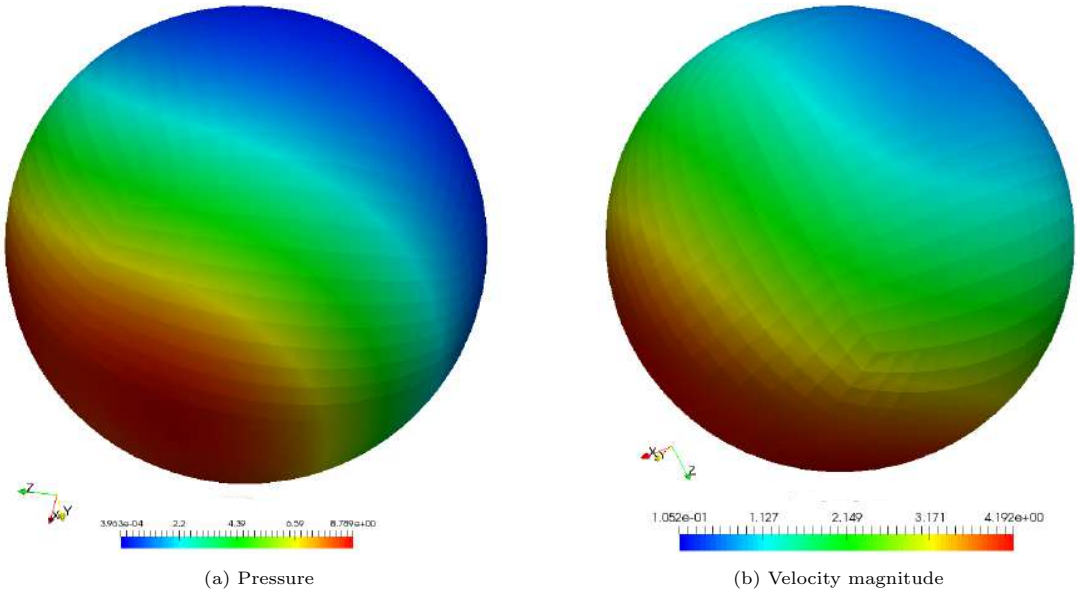

(a) Pressure

(b) Velocity magnitude

Figure 18: Solution of the Ethier-Steinman Navier-Stokes flow in a 3D sphere using $21^3$ quadratic B-spline elements.

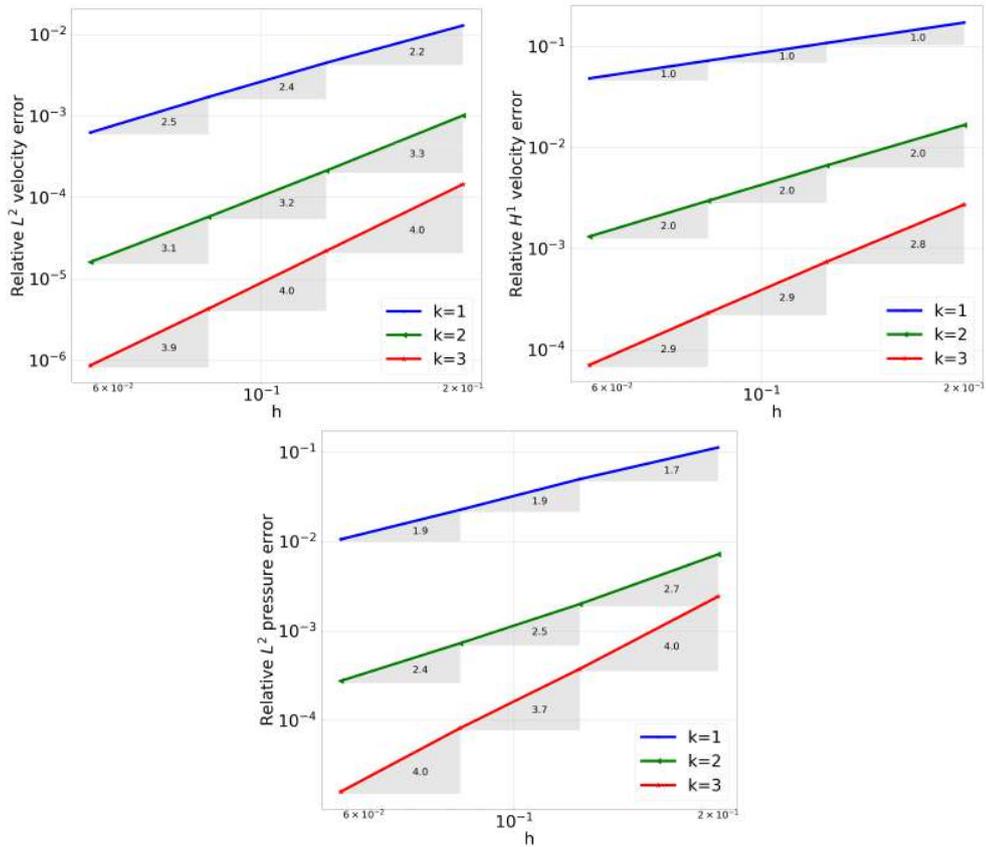

Figure 19: Mesh convergence results for the Ethier-Steinman Navier-Stokes flow in a 3D sphere.

We discretize the problem using a uniform B-spline discretization. In Figure 18 we show the solution obtained using $21^3$ quadratic B-spline elements, from which we observe that the solution is free of oscillations. In Figure 19 we study the mesh convergence behavior for the orders $k = 1, 2, 3$. The considered meshes consist of $5^3$, $8^3$, $12^3$, and $18^3$ elements. We observe optimal rates of converge of $k+1$ and $k$ for the $L^2$-error norm and $H^1$-error norm for the velocity field, respectively. Consistent with the observations of earlier simulations we observe a rate higher than $k$ for the $L^2$-error norm of the pressure field, which we attribute to the use of identical spaces for the pressures and velocities.

## 6. Conclusions

We proposed a stabilization technique for isogeometric analysis of the incompressible Navier-Stokes equations employing the same discretization space for the pressure and velocity fields. The pivotal idea of the developed technique is to penalize the jumps of higher-order derivatives of pressures over element interfaces. Since this technique leverages the skeleton structure of geometric models, we refer to it as a Skeleton-based IsoGeometric Analysis technique. The proposed Skeleton-stabilization penalizes the order $\alpha+1$ derivative jumps for bases with $C^\alpha$ regularity, and hence can be considered as a generalization of continuous interior penalty finite element methods for traditional $C^0$ finite elements. An important advantage of this technique in comparison to inf-sup stable approaches is that it allows the usage of the same discretization space for all field variables. In the context of isogeometric analysis this improves the integration between CAD and analysis, since the technique enables direct usage of the CAD basis for the discretization of all fields.

The proposed Skeleton-Penalty stabilization operator is consistent for solutions with smooth pressure fields. The operator is symmetric and acts only on the pressure space. As a result it does not introduce artificial coupling between the pressure space and the velocity space, and it does not destroy symmetry in the case of the Stokes system. Moreover, no modification of the right-hand-side vector is required, in contrast to some of the alternative stabilization techniques. Considering the bandwidth of the Skeleton-Penalty matrix, there is a substantial advantage to the use of splines, as they ameliorate the large bandwidth that emerges for skeleton-based stabilization operators in Lagrange-based continuous interior penalty methods.

We have observed the proposed Skeleton-Penalty method to yield solutions that are free of pressure oscillations and velocity locking for a wide range of test cases. Optimal convergence rates have been observed for all considered spline orders and regularities, including the case of multi-patch splines. Although a detailed study of the selection of the penalization parameter is beyond the scope of this manuscript, we have observed robustness of the method within a sufficiently large range of penalization parameters.

In this manuscript we have restricted ourselves to the case of moderate Reynolds numbers. Extension to high Reynolds numbers needs a further investigation, as it is anticipated that additional stabilization of the velocity space is then required. We note that in the case of discontinuous spaces – which we have omitted in this work – the proposed stabilization technique fits into the discontinuous Galerkin methodology. We have relied on standard finite element data structures, and we have not considered optimizations that are possible within the isogeometric analysis framework.

**Acknowledgements** We acknowledge the support from the European Commission EACEA Agency, Framework Partnership Agreement 2013-0043 Erasmus Mundus Action 1b, as a part of the *EM Joint Doctorate Simulation in Engineering and Entrepreneurship Development* (SEED). A.R. also acknowledges the support of Fondazione Cariplo - Regione Lombardia through the project "Verso nuovi strumenti di simulazione super veloci ed accurati basati sull'analisi isogeometrica", within the program RST - rafforzamento.

The simulations in this work were performed using the open source software Nutils (www.nutils.org). We would like to thank the Nutils developers for the developments specifically related to this work.